\RequirePackage{fix-cm}
\documentclass[twocolumn]{svjour3}          % twocolumn
\smartqed  % flush right qed marks, e.g. at end of proof

\usepackage{amsmath}
\usepackage{amssymb}
\usepackage{graphicx}
\usepackage{subfigure}
\usepackage[usenames,dvipsnames]{color}
\newcommand{\bv}[1]{\boldsymbol{#1}}
\newcommand{\bx}{\bv{x}}
\newcommand{\bb}{\bv{b}}

\newcommand{\bu}{\bv{u}}
\newcommand{\bp}{\bv{p}}
\newcommand{\bvf}{\bv{f}}
\newcommand{\bg}{\bv{g}}
\newcommand{\bw}{\bv{w}}

\newcommand{\domain}{\Omega}
\newcommand{\bound}{\Gamma}
\newcommand{\R}{\mathbb{R}}
\newcommand{\Rnn}{\R^{N\times N}}
\newcommand{\Rn}{\R^N}
\newcommand{\Rm}{\R^M}
\newcommand{\Rmn}{\R^{M\times N}}
\newcommand{\Rmm}{R^{M\times M}}

\newcommand{\RNum}[1]{\uppercase\expandafter{\romannumeral #1\relax}}

\newcommand{\sect}{Section }
\begin{document}

\title{Fast Multipole Preconditioners for Sparse Matrices Arising from Elliptic Equations%\thanks{Grants or other notes
}


\author{Huda Ibeid \and
        Rio Yokota \and
        Jennifer Pestana \and
        David Keyes
}


\institute{H. Ibeid \at
%             Division of Computer, Electrical and Mathematical Sciences and Engineering\\
           Extreme Computing Research Center\\
           King Abdullah University of Science and Technology,\\ Saudi Arabia\\
           \email{huda.ibeid@kaust.edu.sa}           %  \\
           \and
           R. Yokota \at
%             Division of Computer, Electrical and Mathematical Sciences and Engineering\\
           Tokyo Institute of Technology, Japan\\
           \email{rioyokota@gsic.titech.ac.jp}           %  \\
%%          \emph{Present address:} of F. Author  %  if needed
%%          \emph{Present address:} of F. Author  %  if needed
           \and
           J. Pestana \at
           University of Strathclyde, UK\\
           \email{jennifer.pestana@strath.ac.uk}           %  \\
           \and
           D. Keyes \at
           King Abdullah University of Science and Technology,\\ Saudi Arabia\\
           \email{david.keyes@kaust.edu.sa}           %  \\
}

\date{Received: date / Accepted: date}

\maketitle

\begin{abstract}
Among optimal hierarchical algorithms for the computational solution of elliptic problems, the Fast Multipole Method (FMM) stands out for its adaptability to emerging architectures, having high arithmetic intensity, tunable accuracy, and relaxable global synchronization requirements.  We demonstrate that, beyond its traditional use as a solver in problems for which explicit free-space kernel representations are available, the FMM has applicability as a preconditioner in finite domain elliptic boundary value problems, by equipping it with boundary integral capability for satisfying conditions at finite boundaries and by wrapping it in a Krylov method for extensibility to more general operators. Here, we do not discuss the well developed applications of FMM to implement matrix-vector multiplications within Krylov solvers of boundary element methods. Instead, we propose using FMM for the volume-to-volume contribution of inhomogeneous Poisson-like problems, where the boundary integral is a small part of the overall computation. Our method may be used to precondition sparse matrices arising from finite difference/element discretizations, and can handle a broader range of scientific applications. Compared with multigrid methods, it is capable of comparable algebraic convergence rates down to the truncation error of the discretized PDE, and it offers potentially superior multicore and distributed memory scalability properties on commodity architecture supercomputers. Compared with other methods exploiting the low rank character of off-diagonal blocks of the dense resolvent operator, FMM-preconditioned Krylov iteration may reduce the amount of communication because it is matrix-free and exploits the tree structure of FMM. We describe our tests in reproducible detail with freely available codes and outline directions for further extensibility.
\end{abstract}

\section{Introduction}
\label{intro}
Elliptic PDEs arise in a vast number of applications in scientific computing. A significant class of these involve the Laplace operator, which appears not only in potential calculations but also in, for example, Stokes and Navier-Stokes problems~\cite[Chapters 5 and 7]{elman2005},  electron density computations~\cite[Part II]{martin2004} and reaction-convection-diffusion equations~\cite[Part IV]{hundsdorfer2003}. Conseque-ntly, the rapid solution of PDEs involving the Laplace operator is of wide interest.

Although many successful numerical methods for such PDEs exist, changing computer architectures necessitate new paradigms for computing and the development of new algorithms. Computer architectures of the future will favor algorithms with high concurrency, high data locality, high arithmetic intensity (Flop/Byte), and low synchronicity. This trend is manifested on GPUs and co-processors, where some algorithms are accelerated much less than others on the class of architectures that can be extended to extreme scale. There is always a balance between algorithmic efficiency in a convergence sense, and how well an algorithm scales on parallel architectures. This balance is shifting towards increased parallelism, even at the cost of increasing computation. Since the processor frequency has plateaued for the last decade, Moore's law holds continued promise only for those who are willing to make algorithmic changes.

Among the scientific applications ripe for reconsideration, those governed by elliptic PDEs will be among the most challenging. A common solution strategy for such systems is to discretize the partial differential equations by low-order finite element, finite volume or finite difference methods and then solve the resulting large, sparse linear system. However, elliptic systems are global in nature, and this conflicts with the sweet spots of future architectures. The linear solver must enable the transfer of information from one end of the domain to the other, either through successive local communications (as in many iterative methods), or a direct global communication (as in direct solvers with global recurrences and Krylov methods with global reductions). In either case, avoiding synchronization and reducing communication are the main challenges. There has been considerable effort in this direction in the dense linear algebra community~\cite{demmel2008}. The directed-acyclic-graph-based technology developed in such efforts could be combined with iterative algorithms of optimal complexity for solving elliptic PDEs at extreme scale.

Scalable algorithms for solving elliptic PDEs tend to have a hierarchical structure, as in multigrid methods~\cite{trottenberg2001}, fast multipole methods (FMM)~\cite{greengard1987}, and $\mathcal{H}$-matrices~\cite{hackbusch1999}. This structure is crucial, not only for achieving optimal arithmetic complexity, but also for minimizing data movement. For example, a 3-D FFT requires $\mathcal{O}(\sqrt{P})$ communication for the transpose between pencil-shaped subdomains on $P$ processes \cite{czechowski2012}, whereas these hierarchical methods require $\mathcal{O}(\log P)$ communication \cite{lashuk2012}. This $\mathcal{O}(\log P)$ communication complexity is likely to be optimal for elliptic problems, since an appropriately coarsened representation of a local forcing must somehow arrive at all other parts of the domain for the elliptic equation to converge. In other words, an elliptic problem for which the solution is desired everywhere cannot have a communication complexity of $\mathcal{O}(1)$. However, a disadvantage of these hierarchical methods is that for certain problems we see slow convergence, or even divergence. 

Krylov subspace methods provide another popular alternative to direct methods for general operators. We note that methods such as Chebyshev semi-iteration can require even less communication in the fortunate case when information about the spectrum of the coefficient matrix is known~\cite[Section 10.1.5]{golub1996}, \cite{golub1961}. Among the best known Krylov methods are the conjugate gradient method~\cite{hestenes1952}, MINRES~\cite{paige1975} and GMRES~\cite{saad1986}, although a multitude of Krylov solvers are available in popular scalable solver libraries. The great advantage of these solvers is their robustness -- for any consistent linear system there exists a Krylov method that will converge, in exact arithmetic, for sufficiently many iterations. However, the convergence rate of Krylov methods typically deteriorates as the discretization of an elliptic PDE is refined.

Mesh-independent convergence for Krylov methods applied to systems from elliptic PDEs can often be recovered by preconditioning. Among the best performing preconditioners are the optimal hierarchical methods or, for multiphysics problems such as Stokes and Navier-Stokes equations, block preconditioners with th-ese methods as components. By combining  these hierarchical methods and Krylov subspace solvers we get the benefits of both approaches and obtain a linear solver that is fast but robust. These hierarchical methods all have multiple parameters for controlling the precision of the solution and are able to trade-off accuracy for speed, which is a useful feature for a preconditioner. Furthermore, in analogy to geometric multigrid and algebraic multigrid, $\mathcal{H}^2$-matrices~\cite{hackbusch2000onh2} can be thought of as an algebraic generalization of what FMMs do geometrically. There are advantages and disadvantages to using algebraic and geometric methods, and both have their place as preconditioners. 

There has been some recent work on algebraic multigrids (AMG) in anticipation of the future hardware constraints mentioned above. Gahvari \textit{et al.} developed a performance model for AMG and tested it on various HPC systems -- Intrepid, Jaguar, Hera, Zeus, and Atlas~\cite{gahvari2011}. They found that network distance and contention were both substantial performance bottlenecks for AMG. Adams presents a low-memory matrix-free full multigrid (FMG) with a full approximation storage (FAS)~\cite{adams2012}. He revives an idea from the 1970s~\cite{brandt1977}, which processes the multigrid algorithm vertically, and improves data locality and asynchronicity. Baker \textit{et al.} compared the scalability of different smoothers -- hybrid Gauss-Seidel, $l_1$ Gauss-Seidel, and Chebyshev polynomial, and showed that $l_1$ Gauss-Seidel and Chebychev smoothers scale much better~\cite{baker2012}. There is continuous effort in the multigrid community to adapt the algorithm according to future hardware constraints, and it is likely that multigrid will evolve to remain competitive.

On the other hand, performing a hierarchical low rank approximation (HLRA) of the off-diagonal blocks of a matrix leads to a whole new variety of $\mathcal{O}(N)$ solvers or preconditioners. HLRA based methods include FMM ~\cite{greengard1987}, $\mathcal{H}$-matrices~\cite{hackbusch1999}, hierarchically semi-separable matrices~\cite{chandrasekaran2006}, and recursive skeletonization~\cite{ho2012}. These techniques can be applied to a dense matrix or the Schur complement during a sparse direct solve, thus enabling an $\mathcal{O}(N)$ matrix-vector multiplication of a $N\times N$ dense matrix or an $\mathcal{O}(N)$ direct solve of a $N\times N$ sparse matrix to within a specified accuracy. These HLRA based methods require a smooth kernel in the far field which yields a block low-rank structure. The distinguishing features of the variants come in the way the low rank approximation is constructed -- rank-revealing LU~\cite{pan2000}, rank-revealing QR~\cite{gu1996}, pivoted QR~\cite{kong2011}, truncated SVD~\cite{grasedyck2003}, randomized SVD~\cite{liberty2007}, adaptive cross approximation~\cite{rjasanow2002}, hybrid cross approximation~\cite{borm2005}, and Chebychev interpolation~\cite{dutt1996} are all possibilities. Multipole/local expansions in the FMM constitute another way to construct the low rank approximations. In fact, many of the original developers of FMM are now working on these algebraic variants~\cite{greengard2009}.

Literature on the HLRA based methods mentioned above mainly focuses on the error convergence of the low rank approximation and there is little investigation of the parallel scalability or direct comparison against multigrid. An exception is the work by Grasedyck \textit{et al.}~\cite{grasedyck2008}, where their $\mathcal{H}$-LU preconditioner is compared with BoomerAMG, Pardiso, MUMPS, UMFPACK, SuperLU, and Spooles. However, their executions are serial, and show that their $\mathcal{H}$-matrix code is not yet competitive with these other highly optimized libraries. Another is the work by Gholami \textit{et al.}~\cite{Gholami2014} where they compare FFT, FMM, and multigrids methods for the Poisson problem with constant coefficients on the unit cube with periodic boundary conditions. FMM has also been used as a continuum volume integral with adaptive refinement capabilities~\cite{dhairya2015}. This approach defines the discretization adaptively inside the FMM, whereas in the present method the user defines the discretization and we provide the preconditioner for that given discretization.

In the present work, we consider the Laplace and Stokes equations and devise highly scalable preconditioners for these problems. Our Poisson preconditioner is based on a boundary element method in which matrix-vector multiplies are performed using FMM; the result is an $\mathcal{O}(N)$ preconditioner that is scalable. For the Stokes problem, we apply a block diagonal preconditioner, in which our Poisson preconditioner is combined with a simple diagonal matrix. FMM based preconditioners were first proposed by Sambavaram \textit{et al.} \cite{sambavaram2003}. Such methods lacked practical motivation when flops were expensive, since they turn a sparse matrix into a dense matrix of the same size before hierarchically grouping the off-diagonal blocks. But in a world of cheap flops, the notion of a ``compute-bound preconditioner" sounds more attractive. In the present work, we perform scalability benchmarks and compare the time-to-solution with state-of-the-art multigrid methods such as BoomerAMG in a high performance computing environment.

The rest of the manuscript is organized as follows. In \sect \ref{sec:poisson} we present the model problems and in \sect \ref{sec:krylov} we give an overview of Krylov subspace methods and preconditioning. The basis of our preconditioner is a boundary element method that is discussed in \sect \ref{sec:bem} and the FMM, the essential kernel that makes our method efficient and scalable, is described in \sect    \ref{sec:fmm}. Our numerical results in \sect \ref{sec:results} examine the convergence rates of FMM and multigrid for small Poisson and Stokes problems. Then, in \sect \ref{sec:performance} we scale up the Poisson problem tests and perform strong scalability runs, where we compare the time-to-solution against BoomerAMG \cite{henson2002} on up to 1024 cores. Our conclusions are given in \sect \ref{sec:conc}.

%\begin{figure}
%% Use the relevant command to insert your figure file.
%% For example, with the graphicx package use
%  \includegraphics{example.eps}
%\end{figure}
%\begin{figure*}
%% Use the relevant command to insert your figure file.
%% For example, with the graphicx package use
%  \includegraphics[width=0.75\textwidth]{example.eps}
%\end{figure*}
%\begin{table}
%\begin{tabular}{lll}
%\hline\noalign{\smallskip}
%\end{tabular}
%\end{table}

% Poisson model problem and discretization
\section{Model problems} \label{sec:poisson}
In this section we introduce the Poisson and Stokes model problems we
wish to solve and describe properties of the linear systems that result
from their discretization. We focus on low-order finite elements
but note that discretization by low-order finite difference or finite volume 
methods give linear systems with similar properties.

\subsection{Poisson model problem}
The model Poisson problems we wish to solve are of the form 
\begin{subequations} \label{e:p}
\begin{alignat}{3}
-\nabla\cdot (a \nabla u) & = f && \text{ in } & \domain,\label{e:pint}\\
u & = g && \text{ on } & \bound,\label{e:pbound}
\end{alignat}
\end{subequations}
where $\domain \in \R^d$, $d = 2,3$ is a bounded connected domain with piecewise smooth boundary $\bound$, $f$ is a forcing term, $g$ defines the Dirichlet boundary condition, and $a \geq a_0 > 0 $ is a sufficiently smooth function of space.

Discretization of~\eqref{e:p} by finite elements or finite differences leads to a large, sparse linear system of the form 
\begin{equation} \label{e:axb}
A\bx = \bb,
\end{equation}
where $A\in\Rnn$ is the stiffness matrix and $\bb\in\Rn$ contains the forcing and  boundary data. The matrix $A$ is symmetric positive definite and its eigenvalues  depend on the mesh size, which we denote by $h$, as is typical of discretizations of elliptic PDEs. In particular,  the condition number $\kappa = \lambda_{max}(A)/\lambda_{min}(A)$, the ratio of the largest and smallest eigenvalues of $A$, grows as $O(h^{-2})$ (see, for example,~\cite[Section 1.6]{elman2005}).

\subsection{Stokes model problem}
Incompressible Stokes problems are important when modeling viscous flows and for solving Navier-Stokes equations by operator splitting methods~\cite[Section 2.1]{benzi2005}.  The equations governing the velocity $\bu\in\R^d$, $d=2,3$, and pressure $p\in\R$ of a Stokes fluid in a bounded connected domain $\domain$ with piecewise smooth boundary $\bound$ are \cite{benzi2005},~\cite{elman2005}:
\begin{subequations}\label{e:stokes}
\begin{alignat}{3}
-\nabla^2{\bu} + \nabla p &= {0} && \text{ in } & \domain,\label{e:stokes1} \\ 
\nabla\cdot{\bu} &=0 && \text{ in } & \domain, \label{e:stokes2}\\
{\bu} &= {\bw} && \text{ on } & \bound.\label{e:stokesbc}
\end{alignat}
\end{subequations}

Discretizing~\eqref{e:stokes} by a stabilized\footnote{Although we treat only stabilized discretizations here, stable discretizations are no more difficult to precondition and are discussed in detail in Elman \emph{et al.}~\cite[Chapter 6]{elman2005}.} finite element or finite difference approximation leads to the symmetric saddle point system
\begin{equation}\label{e:sp}
\underbrace{
\begin{bmatrix}
A & B^T\\
B & -C
\end{bmatrix}}_{\mathcal A}
\begin{bmatrix}
\bu\\ \bp
\end{bmatrix}
=
\begin{bmatrix}
\bvf \\ \bg
\end{bmatrix},
\end{equation}
where  ${A}\in\Rnn$ is the vector-Laplacian, a block diagonal matrix with blocks equal to the stiffness matrix from~\eqref{e:axb}, ${B}\in\Rmn$ is the discrete divergence matrix, $C\in\Rmm$ is the symmetric positive definite pressure mass matrix and $\bvf\in\Rn$ and $\bg\in\Rm$ contain the Dirichlet boundary data. 

The matrix $\mathcal A$ is symmetric indefinite
and the presence of the stiffness matrix means that the condition number of $\mathcal{A}$ increases as the mesh is refined. However,  as we will see in the next section, the key ingredient in a preconditioner for $\mathcal A$ that mitigates this mesh dependence is a good preconditioner for the Poisson problem. This allows us to use our preconditioner for the Poisson problem in this more complicated fluid dynamics problem as well. 

% Iterative solvers and preconditioning
\section{Iterative solvers and preconditioning}\label{sec:krylov}
\subsection{Krylov Subspace Methods}
Large, sparse systems of the form~\eqref{e:axb} are often solved by Krylov subspace methods. We focus here on two Krylov methods: the conjugate gradient method (CG)~\cite{hestenes1952} for  systems with symmetric positive definite coefficient matrices and MINRES~\cite{paige1975} for systems with symmetric indefinite matrices. For implementation and convergence details, we refer the reader to the books by Greenbaum~\cite{greenbaum1997} and Saad~\cite{saad2003}.

The convergence of these Krylov subspace methods depends on the spectrum of the coefficient matrix which for the Poisson and Stokes problems, as well as other elliptic PDEs, deteriorates as the mesh is refined. This dependence can be removed by preconditioning. In the case of the Poisson problem~\eqref{e:axb}, we can conceptually think of solving the equivalent linear system $P^{-1}A\boldsymbol{x} = P^{-1}\boldsymbol{b}$ (left preconditioning), or $AP^{-1}\boldsymbol{y} = \bb$, with $P^{-1}\boldsymbol{y} = \bx$ (right preconditioning), for some $P^{-1}\in\Rnn$, and analogously for the Stokes equations~\eqref{e:sp}. However, when the coefficient matrix is symmetric, we would like to preserve this property when preconditioning; this can be achieved by using a symmetric positive definite preconditioner~\cite[Chapters 2 and 6]{elman2005}.  We also note that in practice we never need $P^{-1}$ explicitly but only the action of this matrix on a vector. This enables us to use matrix-free approaches such as multigrid or the fast multipole method.

Many preconditioners for the Poisson problem reduce the number of iterations, with geometric and algebraic multigrid among the most effective strategies~\cite{elman2005}, ~\cite{trottenberg2001}. However, to achieve a lower time-to-solution than can by obtained for the original system, it is also necessary to choose a preconditioner that can be cheaply applied at each iteration. Both geometric and algebraic multigrid methods are $O(N)$, and therefore exhibit good performance on machines and problems for which computation is expensive. However, stresses arise in parallel applications as discussed in the introduction.

For Stokes problems we consider the block diagonal preconditioner
\begin{equation}\label{e:sppre}
{\mathcal P} = 
\begin{bmatrix}
P_A & 0\\
0 & P_S
\end{bmatrix},
\end{equation}
where $P_A\in\Rnn$ and $P_S\in\Rmm$ are symmetric positive definite matrices.
The advantage of this preconditioner is that there is no coupling between the blocks,
so $\mathcal P$ is scalable provided the blocks $P_A$ and $P_S$ are.

Appropriate choices for $P_A$ and $P_S$ have been well studied and it is known that mesh-independent convergence
of MINRES can be recovered when $P_A$ is spectrally equivalent to $A$ in~\eqref{e:sp} and $P_S$ is spectrally equivalent to
the pressure mass matrix $Q\in\Rmm$~\cite{cahouet1988},~\cite[Chapter 6]{elman2005}. These spectral equivalence requirements
imply that the eigenvalues of $P_A^{-1}A$ and $P_S^{-1}Q$ are bounded in an interval on the positive real line independently of the 
mesh width $h$. 

It typically suffices to use the diagonal of $Q$~\cite[Chapter 6]{elman2005}, \cite{wathen1987} or a few steps of Chebyshev semi-iteration ~\cite{wathen2009} for $P_S$. Moreover, the diagonal matrix is extremely parallelizable.
Thus, the key to obtaining a good preconditioner for $\mathcal A$ is to approximate the vector Laplacian effectively. This is typically the most computationally intensive part of the preconditioning process, since in most cases $M \ll N$.

\subsection{The FMM-BEM preconditioner}
In this paper we propose an alternative preconditioner for Poisson and Stokes problems that heavily utilizes the fast multipole method (FMM).  The FMM is $\mathcal{O}(N)$ with compute intensive inner kernels. It has a hierarchical data structure that allows asynchronous communication and execution. These features make the FMM a promising preconditioner for large scale problems on future computer architectures. We show that this preconditioner improves the convergence of these Krylov subspace methods, and is effectively parallelized on today's highly distributed architectures.

The FMM in its original form relies on free-space Green's functions and is able to solve problems with free-field boundary conditions. In Section \ref{sec:bem} the FMM preconditioner is extended to Dirichlet, Neumann or Robin boundary conditions for arbitrary geometries by coupling it with a boundary element method (BEM). Our approach uses the FMM as a \textit{preconditioner} inside a \textit{sparse} matrix solver and the BEM solve is \textit{inside} the preconditioner. Numerous previous studies use FMM for the matrix-vector multiplication inside the Krylov solver for the dense matrix arising from the boundary element discretization. In the present method we are calculating problems with non-zero sources in the volume, and the FMM is used to calculate the volume-to-volume contribution. This means we are performing the action of an $N\times N$ dense matrix-vector multiplication, where $N$ is the number of points in the volume (not the boundary). Additionally, as discussed in Section~\ref{subsec:varcoeff}, it is possible to extend the boundary element method to problems with variable diffusion coefficients, particularly since low accuracy solves are often sufficient in preconditioning.

\begin{figure}
\centering
\includegraphics[width = 0.3\textwidth]{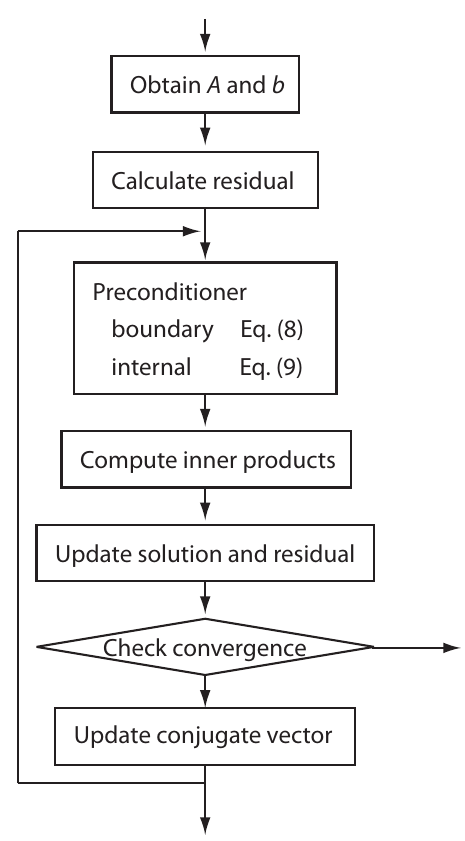}
\caption{Flow chart of the FMM-BEM preconditioner within the conjugate gradient method.}
\label{f:flow_chart}
\end{figure}

Figure \ref{f:flow_chart} shows the flow of calculation of our FMM-BEM preconditioner within the conjugate gradient met-hod; its role in other Krylov solvers is similar. The FMM is used to approximate the matrix-vector multiplication of $A^{-1}$ within the preconditioner. The BEM solver adapts the FMM to finitely applied boundary conditions. During each step of the iteration, the $u$ vector from the previous iteration is used to determine $\partial u/\partial n$ at the boundary from~\eqref{e:ubound}, then~\eqref{e:uinter} is used to compute the new $u$ in the domain $\Omega$.

% BEM preconditioner
\section{Boundary Element Method}\label{sec:bem}
\subsection{Formulation}
We use a standard Galerkin boundary element method \cite{sauter2011} with volume contributions to solve the Poisson equation. A brief description of the formulation is given here. Applying Green's third identity to~\eqref{e:pint} with $a \equiv 1$ gives
\begin{equation}
\int_\Gamma u\frac{\partial G}{\partial n}d\Gamma-\int_\Gamma\frac{\partial u}{\partial n}Gd\Gamma-\int_\Omega u(\nabla^2G)d\Omega=\int_\Omega fGd\Omega,
\label{e:green}
\end{equation}
where $G$ is the Green's function of the Laplace operator, $\frac{\partial}{\partial n}$ is the derivative in the outward normal direction, and $\Gamma$ is the boundary. Following the definition of the Green's function $\nabla^2G=-\delta$, the third term in~\eqref{e:green} becomes
\begin{equation}
-\int_\Omega u(\nabla^2G)d\Omega=\int_\Omega u\delta d\Omega=
\begin{cases}
\frac{1}{2}u\ \text{on}\ \partial\Omega, \\
u\quad \text{in}\ \Omega.
\end{cases}
\end{equation}
Therefore, we may solve the constant coefficient inhomogeneous Poisson problem by solving the following set of equations
\begin{equation}
\int_\Gamma\frac{\partial u}{\partial n}Gd\Gamma=\int_\Gamma u\left(\frac{1}{2}\delta+\frac{\partial G}{\partial n}\right)d\Gamma-\int_\Omega fGd\Omega\quad \text{on}\ \partial\Omega, \label{e:ubound}
\end{equation}

\begin{equation}
u=\int_\Gamma\frac{\partial u}{\partial n}Gd\Gamma-\int_\Gamma u\frac{\partial G}{\partial n}d\Gamma+\int_\Omega fGd\Omega\quad \text{in}\ \Omega. \label{e:uinter}
\end{equation}

As an example, consider the case where Dirichlet boundary conditions are prescribed on $\partial\Omega$. The unknowns are $\partial u/\partial n$ on $\Gamma$ and $u$ in $\Omega\backslash \Gamma$, where~\eqref{e:ubound} solves for the former and~\eqref{e:uinter} can be used to determine the latter. For Neumann boundary conditions one can simply switch the two boundary integral terms in~\eqref{e:ubound} and solve for $u$ instead of $\partial u/\partial n$. In either case, we obtain both $u$ and $\partial u/\partial n$ at each point on the boundary, then calculate~\eqref{e:uinter} to obtain $u$ at the internal points. The last term in~\eqref{e:uinter} takes up most of the calculation time since it is a volume integral for every point in the volume, whereas other terms are either for every point on the boundary or are boundary integrals.

\subsection{Singular integrals}
The Laplace Green's function in 2-D
\begin{equation}
G=-\frac{1}{2\pi}\log r
\label{eq:green}
\end{equation}
is singular. Therefore, the integrals involving $G$ or $\partial G/\partial n$ in~\eqref{e:ubound} and~\eqref{e:uinter} are singular integrals. As described in the following subsection, these singular integral are discretized into piecewise integrals, which are evaluated using Gauss-Legendre quadratures with special treatment for the singular piecewise integral. For boundary integrals in~\eqref{e:ubound} and~\eqref{e:uinter}, analytical expressions exist for the piecewise integral. However, for the volume integral an analytical expression does not exist ~\cite{Aarseth1963}. For this reason, we used a smoothed Green's function of the form
\begin{equation}
G=-\frac{1}{2\pi}\log(\sqrt{r^2+\epsilon^2})
\end{equation}
where $\epsilon$ is a small number that changes with the grid resolution.

\subsection{Discretization}
The integrals in equations~\eqref{e:ubound} and~\eqref{e:uinter} are discretized in a similar fashion to finite element methods. In the following description of the discretization process, we will use the term on the left hand side in~\eqref{e:ubound} as an example. The first step is to break the global integral into a discrete sum of piecewise local integrals over each element
\begin{equation}
\int_\Gamma\frac{\partial u}{\partial n}Gd\Gamma\approx\sum_{j=1}^{N_\Gamma}\int_{\Gamma_j}\frac{\partial u_j}{\partial n}Gd\Gamma_j,
\end{equation}
where $N_\Gamma$ is the number of boundary nodes. These piecewise integrals are performed by using quadratures over the basis functions \cite{sauter2011}. In the present case, we use constant elements so there are no nodal points at the corners of the square domain for the tests in Sections \ref{sec:results} and \ref{sec:performance}. By applying this discretization technique to all terms in~\eqref{e:ubound} we obtain
\begin{multline*}
N_\Gamma
\left\{
\phantom{
\begin{bmatrix}
\ddots\\
&G_{ij}\\
&&\ddots
\end{bmatrix}
}
\right.
\hspace{-24mm}
\overbrace{
\begin{bmatrix}
\ddots\\
&G_{ij}\\
&&\ddots
\end{bmatrix}
}^{N_\Gamma}
\underbrace{
\begin{bmatrix}
\vdots\\
\frac{\partial u_j}{\partial n}\\
\vdots
\end{bmatrix}
}_\text{unknown}
=\\
\overbrace{
\begin{bmatrix}
\hspace{-18mm}\ddots\\
\frac{1}{2}\delta_{ij}+\frac{\partial G_{ij}}{\partial n}\\
\hspace{18mm}\ddots
\end{bmatrix}
}^{N_\Gamma}
\begin{bmatrix}
\vdots\\
u_j\\
\vdots
\end{bmatrix}
-
\overbrace{
\begin{bmatrix}
\ddots\\
&G_{ij}\\
&&\ddots
\end{bmatrix}
}^{N_\Omega}
\begin{bmatrix}
\vdots\\
f_j\\
\vdots
\end{bmatrix},
\end{multline*}
where $N_\Omega$ is the number of internal nodes. All values on the right hand side are known, and $\partial u/\partial n$ at the boundary is determined by solving the linear system. Similarly, we apply the discretization to~\eqref{e:uinter} to have
\begin{multline*}
\small
N_\Omega
\left\{
\begin{bmatrix}
\vdots\\
u_i\\
\vdots
\end{bmatrix}
\right.
=
\overbrace{
\begin{bmatrix}
\hspace{-18mm}\ddots\\
\frac{\partial G_{ij}}{\partial n}\\
\hspace{18mm}\ddots
\end{bmatrix}
}^{N_\Gamma}
\begin{bmatrix}
\vdots\\
u_j\\
\vdots
\end{bmatrix}
\\-
\overbrace{
\begin{bmatrix}
\ddots\\
&G_{ij}\\
&&\ddots
\end{bmatrix}
}^{N_\Gamma}
\begin{bmatrix}
\vdots\\
\frac{\partial u_j}{\partial n}\\
\vdots
\end{bmatrix}
+
\overbrace{
\begin{bmatrix}
\ddots\\
&G_{ij}\\
&&\ddots
\end{bmatrix}
}^{N_\Omega}
\begin{bmatrix}
\vdots\\
f_j\\
\vdots
\end{bmatrix}.
\end{multline*}
At this point, all values on the right hand side are known so one can perform three matrix-vector multiplications to obtain $u$ at the internal nodes, and the solution to the original Poisson equation~\eqref{e:pint}. The third term on the right hand side involves an $N_\Omega\times N_\Omega$ matrix, and is the dominant part of the computational load. This matrix-vector multiplication can be approximated in $\mathcal{O}(N)$ time by using the FMM described in \sect~\ref{sec:fmm}. We also use the FMM to accelerate all other matrix-vector multiplications.

\subsection{Variable coefficient problems}\label{subsec:varcoeff}
A natural question that arises is how to extend the boundary element method, which is the basis of our preconditioner, to problems \eqref{e:p} with variable diffusion coefficients. 

Several strategies for extending boundary element methods to 
problems with variable diffusion coefficients have been proposed 
(see, for example, the thesis of Brunton~\cite[Chapter 3]{brunton1996}). 
Additionally, in this preconditioner setting we may not need to capture the variation
in the diffusion coefficient to a high degree of accuracy; 
for a similar discussion in the context of additive 
Schwarz preconditioners see, for example, Graham \textit{et al.}~\cite{graham2007}. 

Although analytic fundamental solutions can sometimes be found for problems with 
variable diffusion
 (see, e.g., Cheng~\cite{cheng1984} and Clements~\cite{clements1980}), 
in most cases numerical techniques are employed.
One popular method is to introduce a number of subdomains,
on each of which the diffusion coefficient is approximated by a constant 
function~\cite{langer2007,wardle1978}. 

A second option is to split the differential operator into a part for which a fundamental 
solution exists and another which becomes part of the source term. 
Specifically, starting from~\eqref{e:p},
 a similar approach to that described in Banerjee~\cite{banerjee1979} 
and Cheng~\cite{cheng1984} leads to
\begin{multline*}
\int_\bound a u \frac{\partial G}{\partial n}d\bound
 - \int_\Gamma a \frac{\partial u}{\partial n} G 
 - \int_\domain u \nabla a\cdot \nabla G d\domain \\- \int_\domain au\nabla^2G d\domain 
 = \int_\domain fG d\domain,
\end{multline*}

\noindent where again $G$ is the standard fundamental solution for the Laplace operator, i.e., not the fundamental solution for ~\eqref{e:p}.
We can then proceed as described above for~\eqref{e:green}. It is also possible (see Concus and Golub ~\cite{concus}) to change the dependent variable to soak up the variation in $a$ prior to discretization, again resulting in a modified source FEM. 

% FMM kernel
\section{Fast Multipole Method} \label{sec:fmm}
\subsection{Introduction to FMM}
\begin{figure*}
\centering
\subfigure[Direct method]{\includegraphics[width=0.27\textwidth]{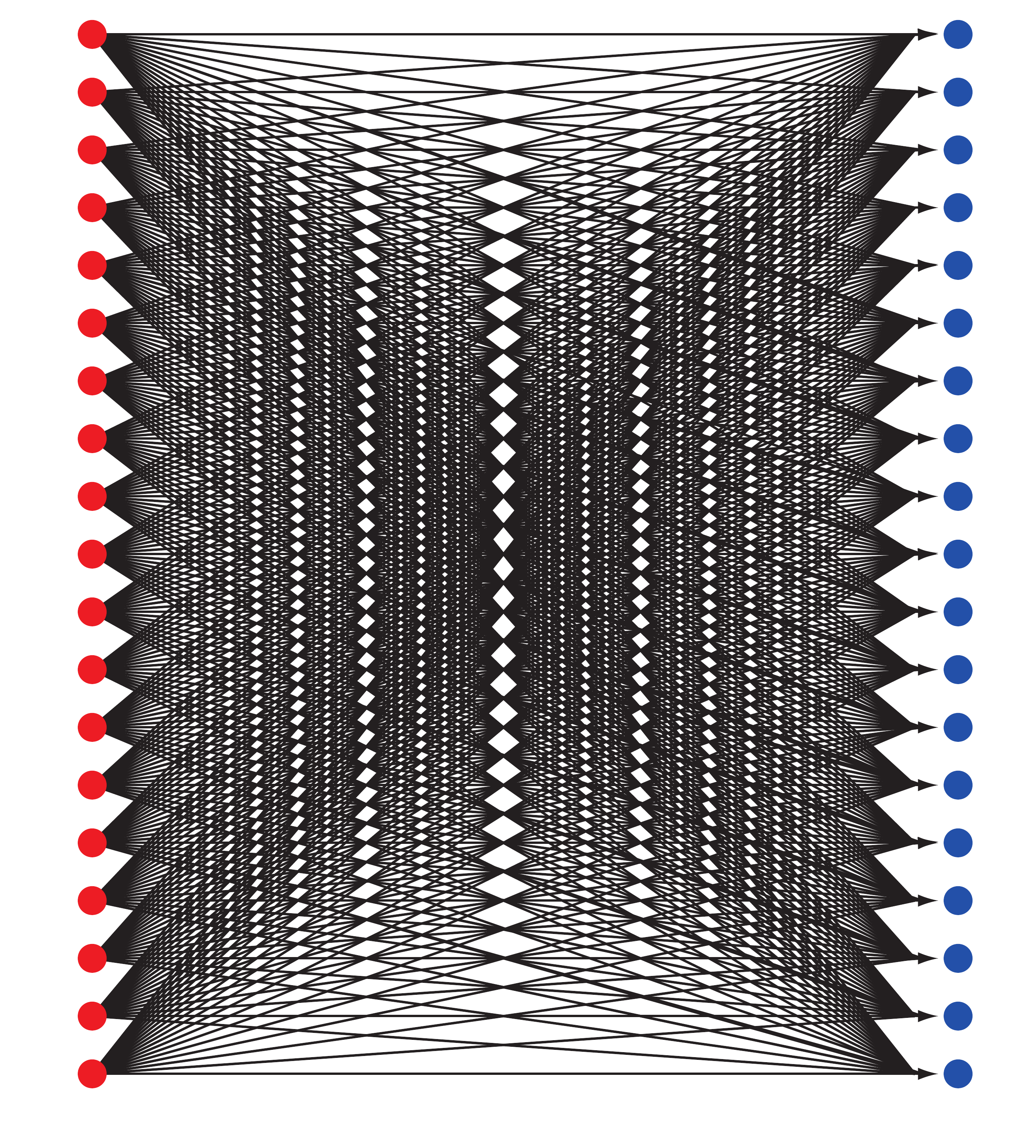}\label{f:direct_interaction}}
\subfigure[Fast Multipole Method]{\includegraphics[width=0.27\textwidth]{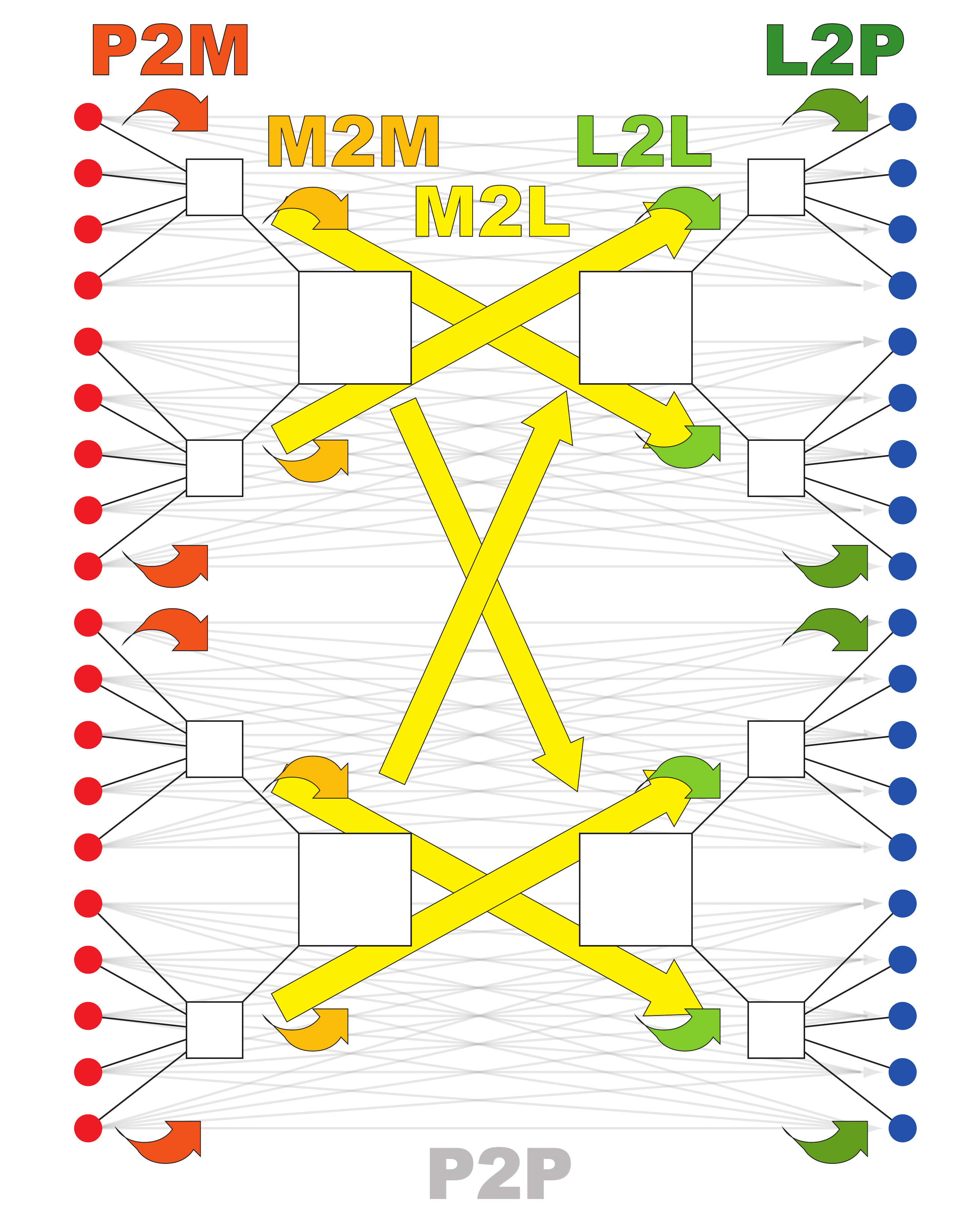}\label{f:fmm_interaction}}
\subfigure[Flow of data in FMM]{\includegraphics[width=0.77\textwidth]{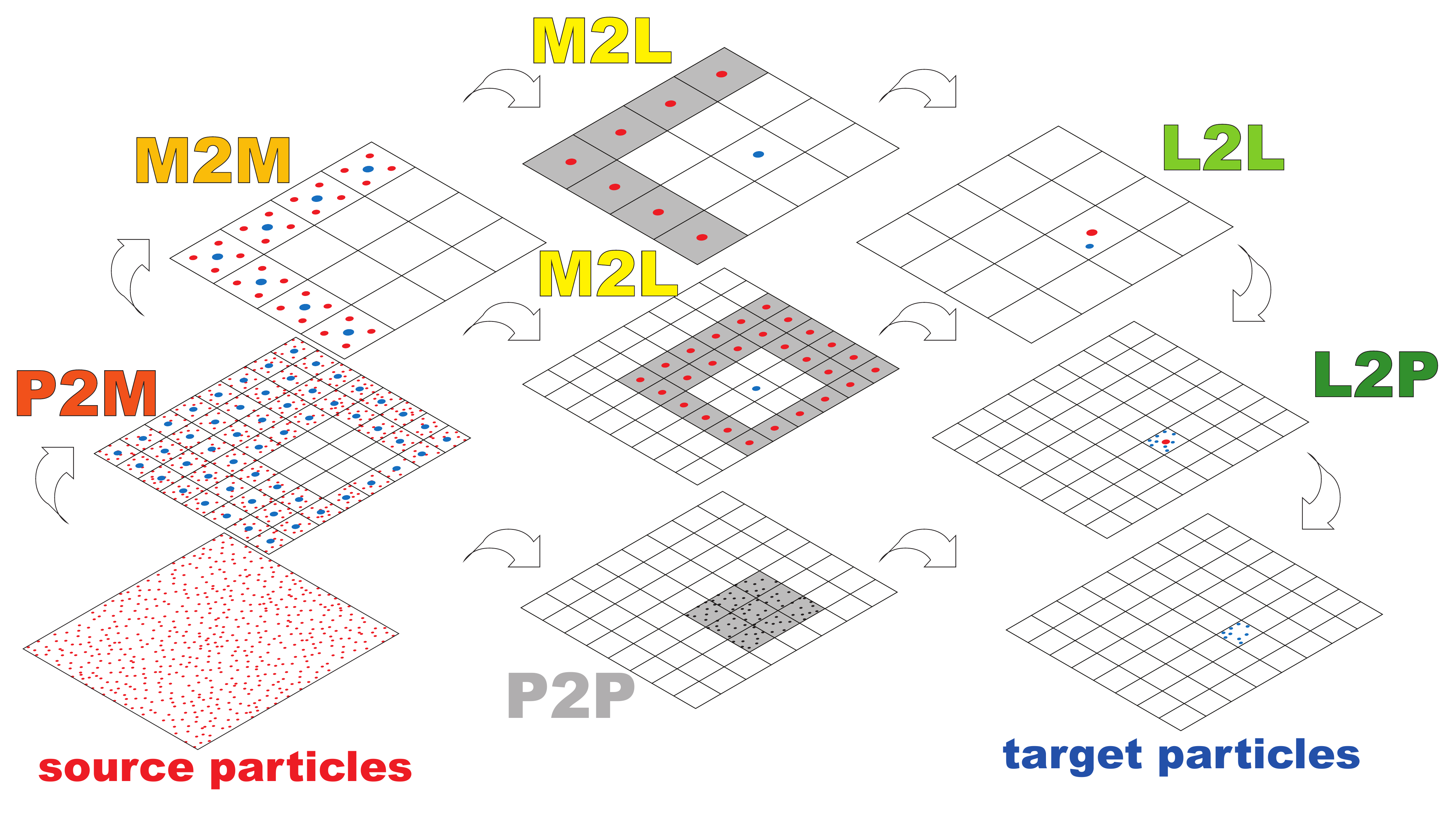}\label{f:fmm_flow}}
\caption{Schematic of Fast Multipole Method. (a) shows the interactions for a $\mathcal{O}(N^2)$ direct method. (b) shows the interactions for the $\mathcal{O}(N)$ FMM, describing the type of interaction between elements in the tree data structure. (c) shows the same FMM kernels as in (b), but from a geometric point of view of the hierarchical domain decomposition.}
\label{f:fmm_schamatic}
\end{figure*}

The last term in Eq.~\eqref{e:uinter} when discretized, has the form
\begin{equation}
u_i=\sum_{j=1}^{N_\Omega}f_jG_{ij}.
\label{e:volume_integral}
\end{equation}
where $i=1,2,...,N_\Omega$. If we calculate this equation directly, it will require $\mathcal{O}(N^2)$ operations. In Figure \ref{f:fmm_schamatic}, we show by schematic how the fast multipole method is able to calculate this in $\mathcal{O}(N)$ operations. Figures \ref{f:direct_interaction} and \ref{f:fmm_interaction} show how the source particles (red) interact with the target particles (blue) for the direct method and FMM, respectively. In the direct method, all source particles interact with all target particles directly. In the FMM, the source particles are first converted to multipole expansions using the P2M (particle to multipole) kernel. Figure \ref{f:fmm_flow} shows the corresponding geometric view of the hierarchical domain decomposition of the particle distribution. Then, multipole expansions are aggregated into larger groups using the M2M (multipole to multipole) kernel. Following this, the multipole expansions are translated to local expansions between well-separated cells using the M2L (multipole to local) kernel. Both Figures \ref{f:fmm_interaction} and \ref{f:fmm_flow} show that the larger cells interact if they are significantly far away, and smaller cells may interact with slightly closer cells. The direct neighbors between the smallest cells are calculated using the P2P (particle to particle) kernel, which is equivalent to the direct method between a selected group of particles. Then, the local expansions of the larger cells are translated to smaller cells using the L2L (local to local) kernel. Finally, the local expansions at the smallest cells are translated into the potential on each particle using the L2P (local to particle) kernel. The mathematical formulae for these kernels will be given in \sect \ref{sec:expansions}. Note, for simplification purposes that each scale of the hierarchical summation can be translated asynchronously from source to target.

In order to perform the FMM calculation mentioned above, one must first decompose the domain in a hierarchical manner. It is common to use an octree in 3-D and quad-tree in 2-D, where the domain is split by its geometrical centerline. The splitting is performed recursively until the number of particles per cell reaches a prescribed threshold. The splitting is usually performed adaptively, so that the densely populated areas result in a deeper branching of the tree. A common requirement in FMMs is that these cells must be isotropic (cubes or squares and not rectangles), since they are used as units for measuring the well-separatedness as shown in Figure \ref{f:fmm_flow} during the M2L interaction. However, our FMM does not use the size of cells to measure the distance between them and allows the cells to be of any shape as long as they can be hierarchically grouped into a tree structure. Once the tree structure is constructed, it is trivial to find parent-child relationships between the cells/particles. This relation is all that is necessary for performing P2M, M2M, L2L, and L2P kernels. However, for the M2L and P2P kernels one must identify a group of well-separated and neighboring cells, respectively. We will describe an efficient method for finding well-separated cells in the following subsection.

\subsection{Dual Tree Traversal}\label{sec:ddt}
The simplest method for finding well-separated pairs of cells in the FMM is to ``loop over all target cells and find their parent's neighbor's children that are non-neighbors," as shown by Greengard and Rokhlin~\cite{greengard1987}. A scheme that permits the interaction of cells at different levels for an adaptive tree was introduced by Carrier \textit{et al.}~\cite{carrier1988}. This scheme is used in many modern FMM codes, and is sometimes called the UVWX-list \cite{lashuk2012}. Another scheme to find well-separated pair of cells is to ``simultaneously traverse the target and source tree while applying a multipole acceptance criterion," as shown by Warren and Salmon~\cite{warren1995}. Teng~\cite{Teng1998} showed that this dual tree traversal can produce interaction pairs that are almost identical to the adaptive interaction list by Carrier \textit{et al.}~\cite{carrier1988}. A concise explanation and optimized implementation of the dual tree traversal is provided by Dehnen~\cite{dehnen2002}.

The dual tree traversal has many favorable properties compared to the explicit construction of interaction lists. First of all, the definition of well-separatedness can be defined quite flexibly. For example, if one were to construct explicit interaction lists by extending the definition of neighbors from $3\times3\times3$ to $5\times5\times5$ using the traditional scheme, the M2L list size will increase rapidly from $6^3-3^3=189$ to $10^3-5^3=875$ in 3-D, which is never faster for any number of expansions. On the other hand, the dual tree traversal can adjust the definition of neighbors much more flexibly and the equivalent interaction list always has a spherical shape. (We say ``equivalent interaction list" because there is no explicit interaction list construction in the dual tree traversal.) The cells no longer need to be cubic, since the cells themselves are not used to measure the proximity of cells. The cells can be any shape or size -- even something like a hierarchical K-means. Of course, the explicit interaction list construction can be modified to include more flexibility, too~\cite{gumerov2008}. However, the resulting code becomes much more complicated than the dual tree traversal, which is literally a few lines of code.\footnote{https://github.com/exafmm/exafmm.git} This simplicity is a large advantage on its own. Furthermore, the parallel version of the dual tree traversal simply traverses the local tree for the target with the local essential tree \cite{warren1992} for the sources, so the serial dual tree traversal code can be used once the local essential tree is assembled.

A possible (but unlikely) limitation of dual tree traversals is the loss of explicit parallelism -- it has no loops. It would not be possible to simply use an OpenMP ``parallel for" directive to parallelize the dual tree traversal. In contrast, the traditional schemes always have an outer loop over the target cells, which can be easily parallelized and dynamically load balanced with OpenMP directives. However, this is not an issue since task based parallelization tools such as Intel Thread Building Blocks (TBB) can be used to parallelize the dual tree traversal. With the help of these tools, tasks are spawned as the tree is traversed and dispatched to idle threads dynamically. In doing so, we not only assure load-balance but also data-locality, so it may actually end up being a superior solution than parallelizing ``for loops" with OpenMP, especially on NUMA architectures.

Considering the advantages mentioned above, we have decided to use the dual tree traversal in our current work. This allows us to perform low accuracy optimizations by adjusting the multipole acceptance criterion without increasing the order of expansions too much, which is the secret to our speed \cite{yokota2013a}. These low accuracy optimizations can give the FMM a performance boost when used as a preconditioner.

\subsection{Multipole Expansions} \label{sec:expansions}
\begin{figure}
\centering
\includegraphics[width = 0.45\textwidth]{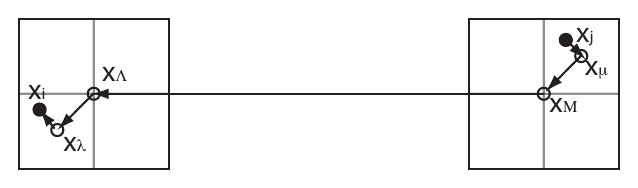}
\caption{Decomposition of the distance vector $\mathbf{x}_{ij}=\mathbf{x}_i-\mathbf{x}_j$ into five parts, that correspond to the five stages P2M, M2M, M2L, L2L, and L2P in the FMM.}
\label{f:vectors}
\end{figure}

For the 2-D Laplace equation, the free space Green's function, as noted in \eqref{eq:green}, has the form
\begin{equation}
G_{ij}=\frac{1}{2\pi}\log\left(\frac{1}{r_{ij}}\right),
\end{equation}
where $r_{ij}=|\mathbf{x}_i-\mathbf{x}_j|$ is the distance between point $i$ and point $j$. By using complex numbers to represent the two-dimensional coordinates $z=x+\iota y$, Eq.~\eqref{e:volume_integral} can be written as
\begin{equation}
u_i=\sum_{j=1}^{N_\Omega}\frac{f_j}{2\pi}\Re\left\{-\log(z_{ij})\right\},
\end{equation}
where $\Re(z)$ represents the real part of $z$. Figure \ref{f:vectors} shows the decomposition of vector $\mathbf{x}_{ij}$ into five parts, $\mathbf{x}_{ij}=\mathbf{x}_{i\lambda}+\mathbf{x}_{\lambda\Lambda}+\mathbf{x}_{\Lambda M}+\mathbf{x}_{M\mu}+\mathbf{x}_{\mu j}$, where $\lambda$ and $\Lambda$ are the center of local expansions and $\mu$ and $M$ are the center of multipole expansions. The lower case is used for the smaller cells and upper case is used for the larger cells. When assuming the relation $|\mathbf{x}_{\Lambda M}|>|\mathbf{x}_{i\lambda}+\mathbf{x}_{\lambda\Lambda}|+|\mathbf{x}_{M\mu}+\mathbf{x}_{\mu j}|$ the following FMM approximations are valid \cite{carrier1988}. We denote the $n$th order multipole expansion coefficient at $\mathbf{x}$ as $M_n(\mathbf{x})$, and the $n$th order local expansion coefficient as $L_n(\mathbf{x})$, where $n=0,1,...,p-1$ for a $p$th order truncation of the series.

\begin{enumerate}
\item P2M from particle at $\mathbf{x}_j$ to multipole expansion at $\mathbf{x}_{\mu}$,
\begin{eqnarray}
M_0(\mathbf{x}_{\mu})&=&\sum_{j=1}^{N}f_j,\\
M_n(\mathbf{x}_{\mu})&=&\sum_{j=1}^{N}\frac{-f_j(-z_{\mu j})^n}{n},  n=\{1,2,...,p-1\}.
\end{eqnarray}

\item M2M from multipole expansion at $\mathbf{x}_\mu$ to multipole expansion at $\mathbf{x}_M$,
\begin{equation}
M_0(\mathbf{x}_M)=M_0(\mathbf{x}_\mu),
\end{equation}
\begin{multline}
M_n(\mathbf{x}_M)=-M_0(\mathbf{x}_\mu)\frac{(-z_{M\mu})^n}{n} \\+\sum_{k=1}^nM_k(\mathbf{x}_\mu)(-z_{M\mu})^{n-k}{{n-1}\choose{k-1}}.
\end{multline}

\item M2L from multipole expansion at $\mathbf{x}_M$ to local expansion at $\mathbf{x}_\Lambda$

\begin{eqnarray}
L_0(\mathbf{x}_\Lambda)\approx M_0(\mathbf{x}_M)\log(z_{\Lambda M})+\sum_{k=1}^{p-1}\frac{M_k(\mathbf{x}_M)}{z_{\Lambda M}^k},
\end{eqnarray}

\begin{multline}
L_n(\mathbf{x}_\Lambda)\approx-\frac{M_0(\mathbf{x}_M)}{(-z_{\Lambda M})^nn}\\+ \sum_{k=1}^{p-1}\frac{(-1)^nM_k(\mathbf{x}_M)}{z_{\Lambda M}^{n+k}}{{n+k-1}\choose{k-1}}.
\end{multline}

\item L2L from local expansion at $\mathbf{x}_\Lambda$ to local expansion at $\mathbf{x}_\lambda$,
\begin{equation}
L_n(\mathbf{x}_\lambda)\approx\sum_{k=n}^{p-1}L_k(\mathbf{x}_\Lambda)z_{\lambda\Lambda}^{k-n}{{k}\choose{n}}.
\end{equation}

\item L2P from local expansion at $\mathbf{x}_\lambda$ to particle at $\mathbf{x}_i$,
\begin{equation}
u_i\approx\Re\left(\sum_{n=0}^{p-1}L_n(\mathbf{x}_{\lambda})z_{i\lambda}^n\right).
\end{equation}
\end{enumerate}
For the P2M, M2M, and M2L kernels, the first term requires special treatment. The expansions are truncated at order $p$, so the accuracy of the FMM can be controlled by adjusting $p$. When recurrence relations are used to calculate the powers of $z$ and the combinations they can be calculated at the cost of one multiplication per inner loop ($k$ loop) iteration. In our implementation, we do not construct any matrices during the calculation of these kernels. The P2P kernel is vectorized with the use of SIMD intrinsics, and the $\log()$ function is calculated using a polynomial fit for $\log_2(x)/(x-1)$ using SIMD.

% Numerical results
\section{Numerical results} \label{sec:results}
In this section we demonstrate the potential of the FMM-based preconditioner by applying it to a number of test problems and comparing it with standard preconditioners. 
The primary aim is to assess the effectiveness of the preconditioner at reducing the number of Krylov subspace iterations that are required for convergence to a given tolerance. Additionally, we seek to ascertain whether mesh independence is achieved. We defer reporting on performance to \sect \ref{sec:performance}. Accordingly, we choose problems that are small enough to enable solution by Matlab.

Our Poisson problems are all two dimensional and include examples with homogeneous and inhomogeneous Dirichlet boundary conditions. We additionally present a two-dimensional Stokes flow problem and show that, as predicted, combining the FMM-based Poisson preconditioner with a block diagonal matrix gives an effective preconditioner for the saddle point problem~\eqref{e:sp}. 

Throughout, our stopping criterion is a decrease in the relative residual of six orders of magnitude. If such a decrease is not achieved after $maxit$  iterations the computations are terminated; this is denoted by  `---' in the tables. This maximum number of iterations is stated for each problem below. For all problems and preconditioners the initial iterate is the zero vector. 

\subsection{The Poisson equation} \label{subsec:poisson}
We first test our preconditioner on three two-dimensional Poisson problems with a constant diffusion coefficient on the domain on $[-1,1]^2$. We discretize the problems by  $Q_1$ finite elements using IFISS~\cite{elman2007},~\cite{ifiss}, with default settings. Our fast multipole preconditioner is compared with the incomplete Cholesky (IC) factorization~\cite{meijerink1977} with zero fill implemented in Matlab and the algebraic multigrid (AMG) and geometric multigrid (GMG) methods in IFISS. Within the GMG preconditioner we select point-damped Jacobi  as a smoother instead of the default ILU, which is less amenable to parallelization. Otherwise, default settings for both multigrid methods are used. For all preconditioners, $maxit=20$ and we apply preconditioned conjugate gradients.

Our first example is the first reference problem in Elman \emph{et al.}~\cite[Section 1.1]{elman2005} for which 
$$
\nabla^2 u = 1 \text{ in } \Omega = [-1,1]^2, \ u = 0 \text{ on } \Gamma.
$$

Table~\ref{t:p1} lists the preconditioned CG iterations for each preconditioner applied. The FMM preconditioner, as well as GMG and AMG appear to give mesh - independent convergence, although the incomplete Cholesky factorization does not.

\begin{table}
\centering
\begin{tabular}{lllll}
\hline\noalign{\smallskip}
$h$ & GMG & AMG & FMM & IC\\
\noalign{\smallskip}\hline\noalign{\smallskip}
$2^{-4}$ &   6 &   5 &   5 &   10 \\
 $2^{-5}$ &   6 &   6 &   6 &   18 \\
 $2^{-6}$ &   7 &   6 &   6 &   --- \\
 $2^{-7}$ &   7 &   6 &   6 &   --- \\
 $2^{-8}$ &   7 &   6 &   7 &   --- \\
 \noalign{\smallskip}\hline
\end{tabular}
\caption{Preconditioned CG iterations for the relative residual to reduce by six orders of magnitude for the problem with $-\nabla^2 u = 1$ and homogeneous boundary conditions.}
\label{t:p1}
\end{table}

\begin{table*}[t]
\centering
\begin{tabular}{lllllll}
\hline\noalign{\smallskip}
$h$ & $\lambda_{min}(A)$ & $\lambda_{max}(A)$ & $\kappa(A)$ & $\lambda_{min}(P^{-1}A)$ & $\lambda_{max}(P^{-1}A)$ & $\kappa(P^{-1}A)$\\
\noalign{\smallskip}\hline\noalign{\smallskip}
$2^{-4}$ & 0.076& 3.94 & 52 & 0.73 & 1.29 & 1.77\\
$2^{-5}$ & 0.019 & 3.99 & 207 & 0.72 & 1.29 & 1.79\\
$2^{-6}$ & 0.005 & 4.00 & 830 & 0.72 & 1.30 & 1.80\\
 \noalign{\smallskip}\hline
\end{tabular}
\caption{Smallest ($\lambda_{min}$) and largest ($\lambda_{max}$) eigenvalues and condition number ($\kappa$) of the stiffness matrix $A$ and FMM-preconditioned matrix $P^{-1}A$ for the problem with $-\nabla^2 u = 1$ and homogeneous boundary conditions.}% with $h = 2^{-4}, \dotsc 2^{-6}$.}
\label{t:p1eig}
\end{table*}

In Table~\ref{t:p1eig} we plot the eigenvalues of the FMM preconditioned stiffness matrix for $h = 2^{-4}, 2^{-5} \text{ and } 2^{-6}$. It is clear that the smallest eigenvalue of $A$ decreases as the mesh is refined; this is particularly problematic for Krylov subspace methods, since small eigenvalues can significantly hamper convergence. However, the eigenvalues of the FMM-preconditioned matrix are bounded away from the origin in a small interval that does not increase in size as the mesh is refined. This hints at spectral equivalence between the FMM-based preconditioner and the stiffness matrix which is unsurprising given that FMM is derived from the exact inverse of the continuous problem. The condition number appears to be bounded, which is a consequence of the mesh-independent convergence observed.

Our second example is the third reference problem from Elman \emph{et al.}~\cite[Section 1.1]{elman2005} posed on $[-1,1]^2$ which is characterized by inhomogeneous Dirichlet boundary conditions and the analytic solution 
$$u(x,y) = \frac{2(1+y)}{(3+x)^2 + (1+y)^2}.$$  
From Table~\ref{t:p2} we find that,  similarly to the previous problem, the FMM preconditioner and both multigrid preconditioners are mesh independent but the Cholesky preconditioner is not. The FMM preconditioner is also competitive with the multigrid methods. Thus, on systems on which applying the FMM preconditioner is significantly faster than applying the multigrid preconditioners, we will achieve a faster time-to-solution with the former. We note that the eigenvalues and condition numbers obtained for the FMM preconditioned stiffness matrix are the same as those computed for the previous example. 

\begin{table}
\centering
\begin{tabular}{lllll}
\hline\noalign{\smallskip}
$h$ & GMG & AMG & FMM & IC\\
\noalign{\smallskip}\hline\noalign{\smallskip}
$2^{-4}$ &   5 &   5 &   5 &   11 \\
 $2^{-5}$ &   5 &   5 &   5 &   19 \\
 $2^{-6}$ &   5 &   5 &   5 &   --- \\
 $2^{-7}$ &   5 &   5 &   5 &   --- \\
 $2^{-8}$ &   5 &   5 &   5 &   --- \\
 \noalign{\smallskip}\hline
\end{tabular}
\caption{Preconditioned CG iterations for the relative residual to reduce by six orders of magnitude for the problem with $-\nabla^2 u = 0$ and inhomogeneous boundary conditions.}
\label{t:p2}
\end{table}

The final problem we consider in this section is the Poisson problem with solution
$$u(x,y) = x^2 + y^2$$
on $[-1,1]^2$, which has forcing term $f\equiv -4$ in the domain and inhomogeneous Dirichlet boundary conditions. The convergence results for this problem, given in Table~\ref{t:p3}, are similar to those for the previous problems. They show that the FMM preconditioner gives mesh independent convergence and is competitive with AMG and GMG. We also obtain the same eigenvalue results as for the previous  examples.

\begin{table}
\centering
\begin{tabular}{lllll}
\hline\noalign{\smallskip}
$h$ & GMG & AMG & FMM & IC\\
\noalign{\smallskip}\hline\noalign{\smallskip}
$2^{-4}$ &   5 &   5 &   5 &   10 \\
 $2^{-5}$ &   5 &   5 &   5 &   18 \\
 $2^{-6}$ &   5 &   5 &   5 &   --- \\
 $2^{-7}$ &   5 &   5 &   5 &   --- \\
 $2^{-8}$ &   5 &   5 &   5 &   --- \\
 \noalign{\smallskip}\hline
\end{tabular}
\caption{Preconditioned CG iterations for the relative residual to reduce by six orders of magnitude for the problem with $-\nabla^2 u = -4$ and inhomogeneous  boundary conditions.}
\label{t:p3}
\end{table}

\begin{figure}
\centering
\includegraphics[width = 0.47\textwidth,height=7cm]{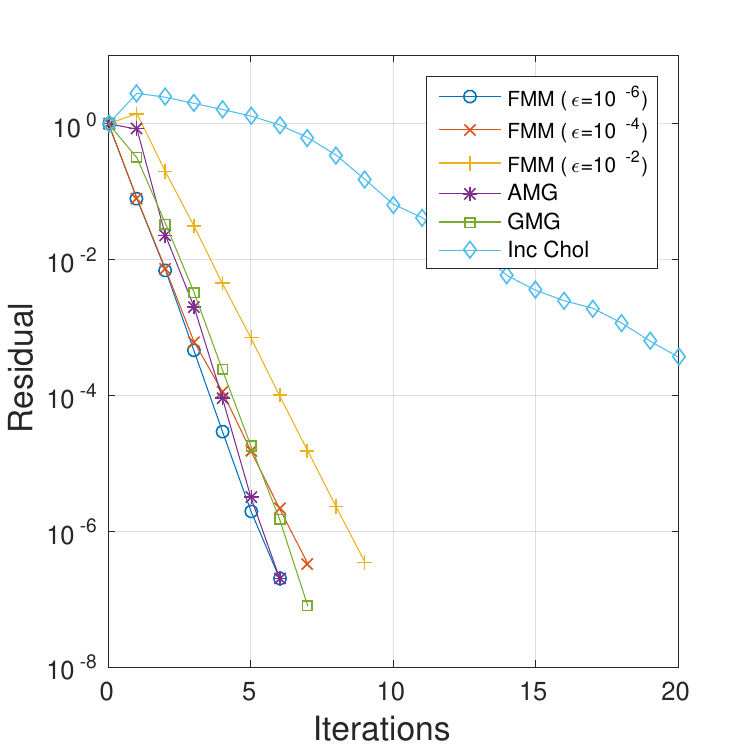}
\caption{Convergence rate of the FMM preconditioner with different precision, plotted along with algebraic multigrid, geometric multigrid, and incomplete Cholesky preconditioners. The $\epsilon$ represents the precision of the FMM, where $\epsilon=10^{-6}$ corresponds to six significant digits of accuracy.}
\label{f:convergence}
\end{figure}

\subsection{Effect of FMM precision on convergence}
For the results shown above, the FMM precision was set to preserve six significant digits. However, the FMM can be accelerated further by trading precision for speed. Since we are using the FMM as a preconditioner, the accuracy requirements are somewhat lower than that of general applications of FMM. Although this balance between the accuracy and speed of FMM is a critical factor for evaluating the usefulness of FMM as a preconditioner, the relation between the FMM precision and convergence rate has not been studied previously.

In Figure \ref{f:convergence} the relative residual at each CG iteration is plotted against the number of iterations for FMM, AMG, GMG, and IC. The problem is the same as in Table \ref{t:p1}. Three cases of FMM are used with six, four, and two significant digits of accuracy, respectively. The $\epsilon=10^{-6}$ case corresponds to the condition for the tests in Tables \ref{t:p1}--\ref{t:p3}. Decreasing the FMM accuracy to four digits has little effect during the first few iterations, but slows down the convergence near the end. Decreasing the FMM accuracy further to two digits slows down the convergence somewhat, but is still much better than the incomplete Cholesky.

\begin{figure}
\centering
\includegraphics[width = 0.5\textwidth]{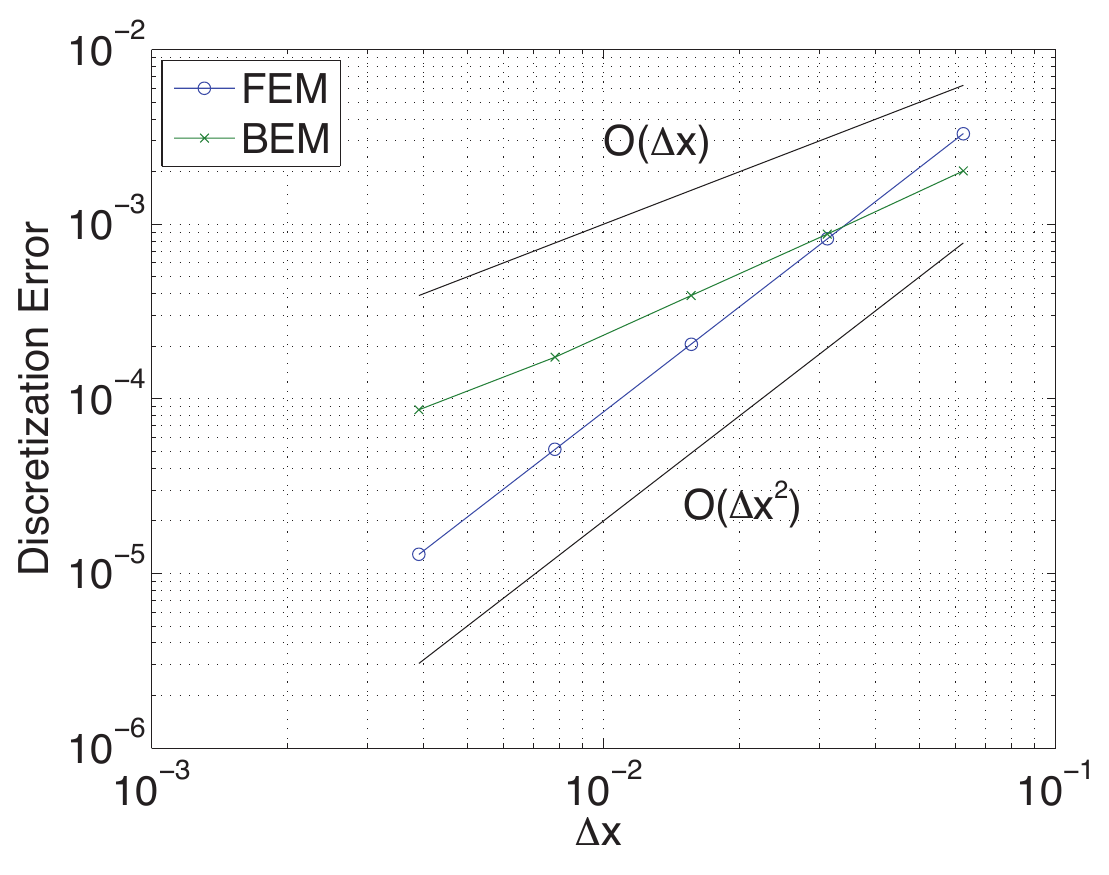}
\caption{Convergence of spatial discretization error for the FEM and BEM. The relative $L^2$ norm of the difference between the analytical solution is plotted against the grid spacing $\Delta x$.}
\label{f:discretization}
\end{figure}

Increasing the precision of the FMM past six digits did not result in any noticeable improvement because truncation error begins to dominate. We are preconditioning a matrix resulting from a FEM discretization by using a integral equation with Green's function kernels. Each has its own error, below which algebraic error need not be reduced. We show in Figure \ref{f:discretization} the convergence of spatial discretization error for the FEM and BEM approaches. We use the same reference problem as in Table \ref{t:p1}, which has an analytical solution. The discretization error is measured by taking the relative $L^2$ norm of the difference between the analytical solution and the individual numerical solutions. We see that the FEM is second order and BEM is first order. The five different values of $\Delta x$ correspond to $h=\{2^{-4},2^{-5},2^{-6},2^{-7},2^{-8}\}$, which were used in the previous experiments. For the current range of grid spacing, the discrepancies between the FEM and BEM truncation error is in the range of $10^{-3}$ to $10^{-4}$.

\subsection{Stokes problem}
Next, we examine convergence for a two-dimensional Stokes flow. The leaky cavity problem~\cite[Example 5.1.3]{elman2005} on $[-1,1]$ is discretized by $Q_1-P_0$ elements in Matlab using IFISS with default settings. As described in Section~\ref{sec:krylov}, by combining a stiffness matrix preconditioner $P_A$ with the diagonal of the pressure mass matrix $P_S$, an effective preconditioner~\eqref{e:sppre} for the saddle point system~\eqref{e:sp} is obtained. Here, we are interested in using the FMM preconditioner for $P_A$, and we compare its performance with AMG and GMG. We do not consider the incomplete Cholesky factorization of $A$ because of its poor performance on the stiffness matrix (see Tables~\ref{t:p1}, \ref{t:p2} and~\ref{t:p3}).  We set $maxit=50$ and apply preconditioned MINRES to the saddle point system.

As for the Poisson problem, the FMM-based preconditioner provides a mesh-independent preconditioner that is comparable to algebraic and geometric multigrid. Although two or three more iterations are required by the FMM preconditioner than the AMG preconditioner, if each iteration is faster the time-to-solution may be lower.

\begin{table}
\centering
\begin{tabular}{llll}
\hline\noalign{\smallskip}
$h$ & GMG & AMG & FMM\\
\noalign{\smallskip}\hline\noalign{\smallskip}
$2^{-4}$ &   32 &   31 &   35 \\
 $2^{-5}$ &   32 &   32 &   35 \\
 $2^{-6}$ &   33 &   32 &   33 \\
 $2^{-7}$ &   31 &   31 &   33 \\
 $2^{-8}$ &   31 &   31 &  31\\
\noalign{\smallskip}\hline
\end{tabular}
\caption{Preconditioned MINRES iterations for the relative residual to reduce by six orders of magnitude for the Stokes problem.}
\label{t:s1}
\end{table}

\subsection{Variable coefficient Poisson equation}

For sufficiently smooth diffusion coefficient variation, we can precondition the variable coefficient problem with the constant coefficient problem, since they are spectrally equivalent. We test this approach on the variable coefficient Poisson equation of the form 
\begin{subequations}
\begin{alignat*}{3}
-\nabla \cdot (a \nabla u) & = 1 && \text{ in } & \domain,\\
u & = 0 && \text{ on } & \bound,
\end{alignat*}
\end{subequations}
where
\begin{equation*}
a = 1+\mu(\sin(m \pi x)\sin(n \pi y)).
\end{equation*}

Tables~\ref{t:v1},~\ref{t:v2}, and~\ref{t:vv2}  show that the FMM preconditioner and both multigrid preconditioners achieved mesh independent convergence for different amplitudes $\mu$ and frequencies $m$ and $n$. Also, the FMM preconditioner is competitive with the algebraic multigrid method requiring comparable number of iterations.

Similar to Table~\ref{t:p1eig}, Table \ref{t:v3} shows that the eigenvalues of the FMM-preconditioned matrix are bounded away from the origin.

\begin{table}
\centering
\begin{tabular}{llll}
\hline\noalign{\smallskip}
$\mu$ & AMG & FMM\\
\noalign{\smallskip}\hline\noalign{\smallskip}
 $2^{-16}$ &   6 &   6 \\
 $2^{-8}$ &   6 &   6 \\
 $2^{-4}$ &   6 &   7 \\   
 $2^{-2}$ &   6 &   8 \\   
 \noalign{\smallskip}\hline
\end{tabular}
\caption{Preconditioned CG iterations for the relative residual to reduce by six orders of magnitude ($h = 2^{-6}$, $m = 1$ and $n = 1$).}
\label{t:v1}
\end{table}

\begin{table}
\centering
\begin{tabular}{llll}
\hline\noalign{\smallskip}
$m$ & $n$ &  AMG & FMM\\
\noalign{\smallskip}\hline\noalign{\smallskip}
1 & 1 &   6 &   7 \\
2 & 2 &   6 &   7 \\
4 & 4 &   6 &   7 \\
8 & 8 &   6 &   7 \\   
16 & 16 &   6 &   7 \\     
 \noalign{\smallskip}\hline
\end{tabular}
\caption{Preconditioned CG iterations for the relative residual to reduce by six orders of magnitude ($h = 2^{-6}$, $\mu = 2^{-4}$).}
\label{t:v2}
\end{table}

\begin{table}
\centering
\begin{tabular}{lll}
\hline\noalign{\smallskip}
$n$ &  AMG & FMM\\
\noalign{\smallskip}\hline\noalign{\smallskip}
1 &   6 &   7 \\
2 &   6 &   7 \\
4 &   6 &   7 \\
8 &   6 &   7 \\   
16 &   6 &   7 \\     
 \noalign{\smallskip}\hline
\end{tabular}
\caption{Preconditioned CG iterations for the relative residual to reduce by six orders of magnitude ($m = 4$, $h = 2^{-6}$, $\mu = 2^{-4}$).}
\label{t:vv2}
\end{table}

 \begin{table*}[t]
\centering
\begin{tabular}{llllllllll}
\hline\noalign{\smallskip}
$h$ & $m$ & $n$ & $\lambda_{min}(A)$ & $\lambda_{max}(A)$ & $\kappa(A)$ & $\lambda_{min}(P^{-1}A)$ & $\lambda_{max}(P^{-1}A)$ & $\kappa(P^{-1}A)$ \\
\noalign{\smallskip}\hline\noalign{\smallskip}
$2^{-4}$ & 3 & 3 &   0.0759 & 3.9923 & 53 &   0.4423 &  1.0000 &  2.26\\
$2^{-5}$ & 6 & 6 &   0.0192 & 4.0199 & 209 &   0.4371 &  1.0040 &  2.29\\
$2^{-6}$ & 12 & 12 &  0.0048 & 4.0280 & 839 &   0.4360 & 1.0061 & 2.31\\
 \noalign{\smallskip}\hline
\end{tabular}
\caption{Smallest ($\lambda_{min}$) and largest ($\lambda_{max}$) eigenvalues and condition number ($\kappa$) of the stiffness matrix $A$ and FMM-preconditioned matrix $P^{-1}A$ with $\mu = 2^{-4}$.}
\label{t:v3}
\end{table*}

% \begin{table}[h]
%\centering
%\begin{tabular}{crrrrrrrrr}
%\hline
%$h$ & $m$ & $n$ & $\lambda_{min}(A)$ & $\lambda_{max}(A)$ & $\kappa(A)$ & $\lambda_{min}(P^{-1}A)$ & $\lambda_{max}(P^{-1}A)$ & $\kappa(P^{-1}A)$ \\
%\hline
%$2^{-6}$ & 3 & 3 &  0.0048 & 4.1690 & 869 &   0.4209 &  1.0338 &  2.46\\
%$2^{-7}$ & 6 & 6 &  0.0012 &   4.1691 &   3474 &  0.4209 &  1.0338 &  2.46\\
%$2^{-8}$ & 12 & 12 &  0.0003 &   4.1692 &   13427 &  ---&  ---&  ---\\
% \hline
%\end{tabular}
%\caption{Smallest ($\lambda_{min}$) and largest ($\lambda_{max}$) eigenvalues and condition number ($\kappa$) of the stiffness matrix $A$ and FMM-preconditioned matrix $P^{-1}A$ with $\mu = 2^{-4}$.}
%\label{t:v3}
%\end{table}

\section{Performance analysis}\label{sec:performance}
In this section we evaluate the performance of the FMM-based preconditioner by comparing its time-to-solution to an algebraic multigrid code BoomerAMG. We have implemented our FMM-preconditioner in PETSc~\cite{petsc-user-ref}~,\cite{petsc-web-page} via PetIGA~\cite{PetIGA}. PetIGA is a software layer that sits on top of PETSc that facilitates NURBS-based Galerkin finite element analysis. For our present analysis, we simply use PetIGA to reproduce the same finite element discretization as the tests shown in Section \ref{sec:results}, but in a high performance computing environment. We select the first problem in Section \ref{subsec:poisson} with $-\nabla^2u=1$ and homogeneous Dirichlet boundary conditions for the following performance evaluation.

All codes that were used for the current study are publicly available. A branch of PetIGA that includes the FMM preconditioner is hosted on bitbucket. \footnote{https://bitbucket.org/rioyokota/petiga-fmm}
% The 2-D extension of exaFMM, which we call from PetIGA is also hosted on bitbucket. \footnote{https://bitbucket.org/rioyokota/exafmm2d}

All calculations were performed on the TACC Stampede system without using the coprocessors. Stampede has 6400 nodes, each with 
 two Xeon E5-2680 processors and one Intel Xeon Phi SE10P coprocessor and 32GB of memory. We used the Intel compiler (version 13.1.0.146) and configured PETSc with ``\texttt{COPTFLAGS=-O3 FOPTFLAGS=-O3 --with-clanguage=cxx\\ --download-f-blas-lapack --download-hypre\\ --download-metis --download-parmetis\\ --download-superlu\_dist --with-debugging=0}".

\begin{figure}[t]
\centering
\includegraphics[width=0.5\textwidth,height=6cm]{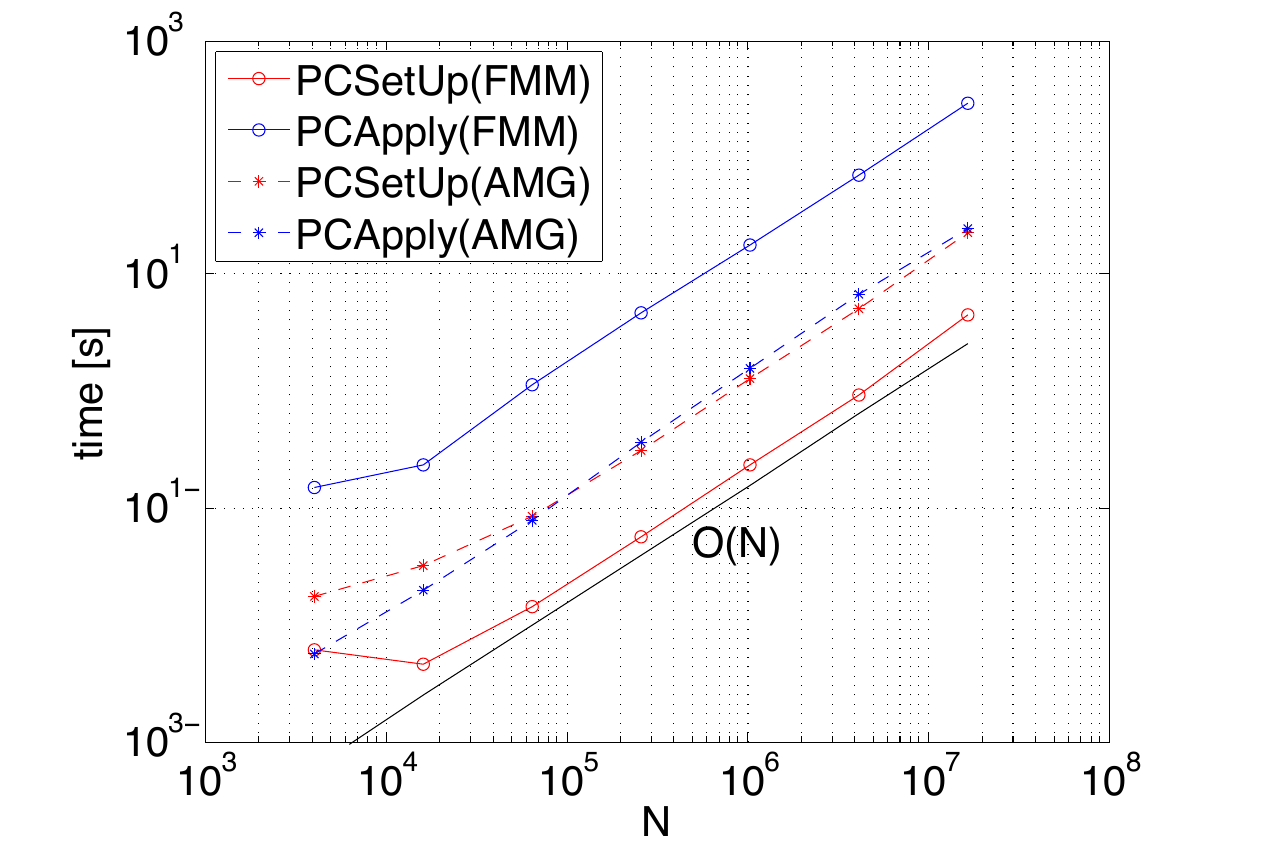}
\caption{Time-to-solution for different problem sizes of the FMM and AMG preconditioners on a single core of a Xeon E5-2680.}
\label{fig:complexity}
\end{figure}

\subsection{Serial results}
We first evaluate the serial performance of our method using the same two-dimensional Poisson problem used in Section \ref{sec:poisson}. We confirmed that the iteration counts shown in Table \ref{t:p1} did not change for the PETSc version of our code. Then, we measured the time-to-solution for different problem sizes. Since the domain size is $[-1,1]$, the grid spacing of $h=\{2^{-4},2^{-5},2^{-6},2^{-7},2^{-8}\}$ in Section \ref{sec:results} correspond to a grid size of $N=\{32^2,64^2,128^2,\\ 256^2,512^2\}$. In the PETSc version, the time-to-solution improves significantly so we tested for larger problem sizes of $N=\{64^2,128^2,256^2,512^2,1024^2,2048^2,4096^2\}$.

The time-to-solution is plotted against the problem size $N$ in Figure \ref{fig:complexity}. Since we are using PETSc, it is trivial to change the preconditioner to AMG by passing the option ``\texttt{--pc\_type hypre}" during runtime. Therefore, the time-to-solution of BoomerAMG is shown as a reference in the same figure. For BoomerAMG we compared different relaxation, coarsening, and interpolation methods and found that \\ ``\texttt{-pc\_hypre\_boomeramg\_relax\_type\_all\\ backward-SOR/Jacobi\\ -pc\_hypre\_boomeramg\_coarsen\_type\\ modifiedRuge-Stueben\\ -pc\_hypre\_bommeramg\_interp\_type classical}" gives the best performance.

Both FMM and AMG runs are serial, where we used a single MPI process and a single thread. The majority of the time goes into the setup of the preconditioner ``PCSetUp" and the actual preconditioning ``PCApply", so only these events are shown in the legend. The ``PCSetUp" is called only once for the entire run, while ``PCApply" is called every iteration. For the present runs, both FMM and AMG required six iterations for the relative residual to drop six digits, so all runs are calling ``PCApply" six times. The order of expansion for the FMM is set to $p=6$ and $\theta=0.4$, which gives about six significant digits of accuracy. With this accuracy for the FMM, we are still able to converge in six iterations. The P2P kernel in the FMM code is performed in single precision using SIMD intrinsics, but this does not prevent us from reaching the required accuracy of six significant digits because we use Kahan's summation technique \cite{kahan1965} for the reduction.

By taking a closer look at Figure \ref{fig:complexity}, one can see that both the FMM and AMG show $\mathcal{O}(N)$ asymptotic behavior. The FMM seems to have a slower preconditioning time, but a much faster setup time compared to AMG. The FMM also has a constant overhead which becomes evident when $N$ is small. In summary, the time-to-solution of the FMM is approximately an order of magnitude larger than that of AMG for the serial runs. This is consistent with our intuition that FMM is not the preconditioner of choice for solving small problems on a single core. We will show in the following section that the FMM becomes competitive when scalability comes into the picture.

\begin{figure}[t]
\centering
\includegraphics[width=0.48\textwidth,height=7cm]{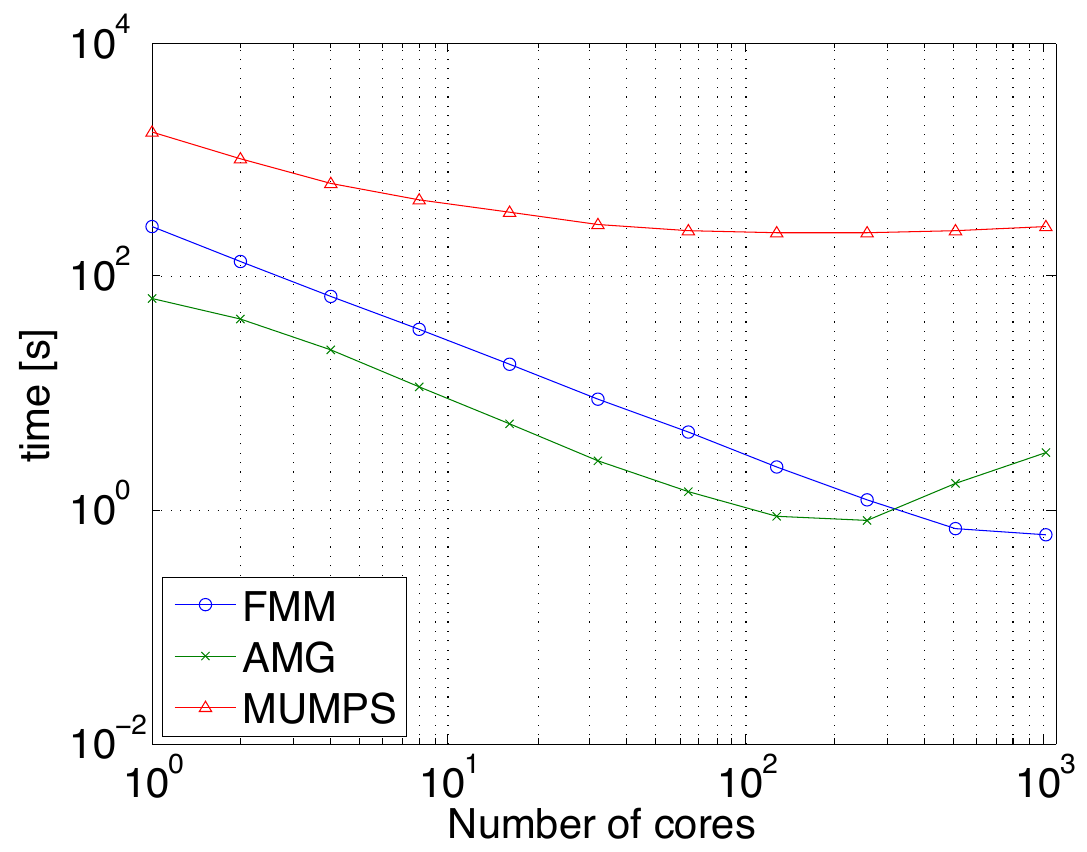}
\caption{Strong scaling of the 2-D FMM and AMG preconditioners.}
\label{fig:strong_scaling_2d}
\end{figure}

\subsection{Parallel results}
Using the same Poisson problem, we now compare the performance of FMM and AMG for parallel runs on Stampede. We also compare with a sparse direct solver MUMPS by invoking at runtime ``\texttt{-ksp\_type\ preonly -pc\_type lu -pc\_factor\_mat\_solver\_package mumps}".

The strong scaling of FMM, AMG, and MUMPS are shown in Figure \ref{fig:strong_scaling_2d}. We use the largest grid size in the previous runs $N=4096^2$. Stampede has 16 cores per node so all runs first parallelize among these cores and then among the nodes after the 16 cores are filled. The FMM strong scales quite well up to 1024 cores, while the parallel efficiency of AMG starts to decrease after 128 cores. The sparse direct solver has a much larger time-to-solution even on a single core, and is much less scalable than the other two hierarchical preconditioners. It is worth mentioning that the setup cost of the direct solver is dominant and so if several linear systems are solved with the same coefficient matrix then this cost is amortized. For this particular Poisson problem on this particular machine using this particular FMM code we see an advantage over BoomerAMG past 512 cores.

\subsection{Extension to 3-D}
The results above are all two-dimensional. A natural question that arises is whether the extension to 3-D is straightforward, and whether FMM will still be competitive as a preconditioner or not. Our results showed that a dominant part of the calculation time for the FMM preconditioner is the ``PCApply" stage, which is the dual tree traversal for calculation of M2L and P2P kernels. For 3-D kernels, the M2L operation is much more complicated so the calculation time of the FMM will increase, even for the same number of unknowns $N$.

\begin{figure}[t]
\centering
\includegraphics[width=0.45\textwidth,height=7cm]{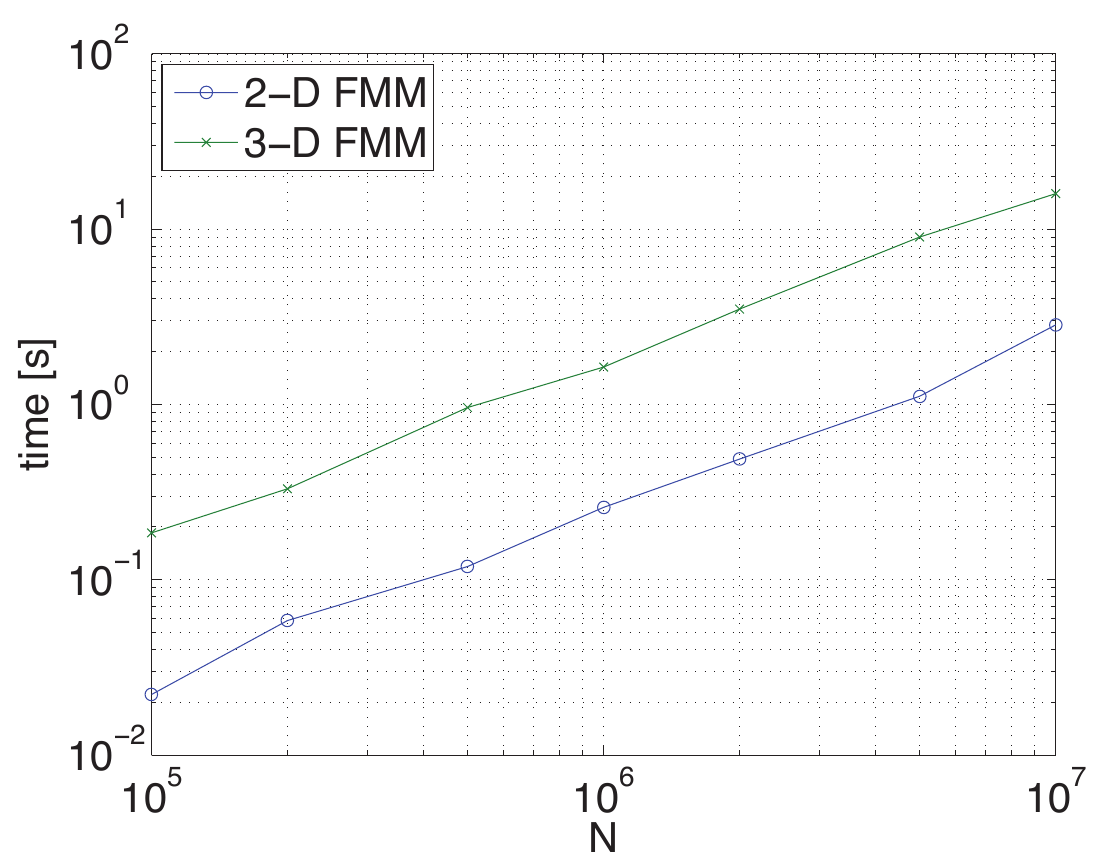}
\caption{Calculation time of 2-D and 3-D FMM for the same problem size.}
\label{f:2dvs3d}
\end{figure}

\begin{figure}[t]
\centering
\includegraphics[width=0.5\textwidth,height=7cm]{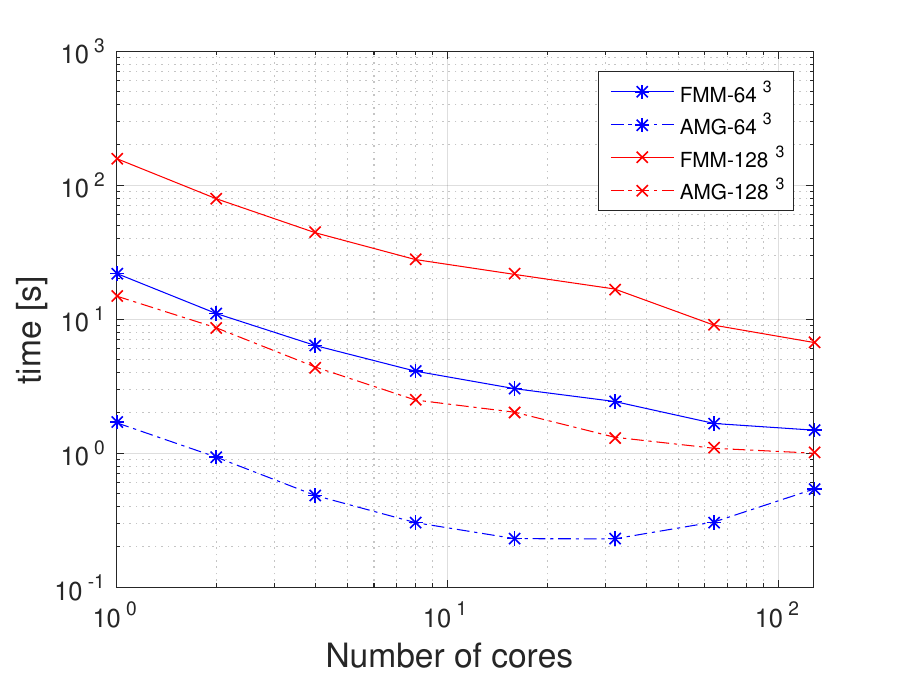}
\caption{Strong scaling of the 3-D FMM and AMG preconditioners.}
\label{fig:strong_scaling_3d}
\end{figure}

Figure \ref{f:2dvs3d} shows the calculation time of our 2-D FMM and 3-D FMM, both for the Laplace kernel with four significant digits of accuracy on a single core of a Xeon E5-2680, 2.7 GHz CPU. The problem size $N$ varies from $10^5$ to $10^7$. We see that the 3-D FMM is about an order of magnitude slower than the 2-D FMM for the same problem size. Nevertheless, Figure~\ref{fig:strong_scaling_3d} shows that the 3-D FMM preconditioner strong scales quite well up to 128 cores for $N=64^3$ and $N=128^3$ when compared to BoomerAMG with these configurations\\``\texttt{-pc\_hypre\_boomeramg\_coarsen\_type hmis\\ -pc\_hypre\_boomeramg\_interp\_type ext+i\\ -pc\_hypre\_boomeramg\_p\_max 4\\ -pc\_hypre\_boomeramg\_agg\_nl 1}". \\
These runs were performed on Shaheen \RNum{2} which is a Cray XC40 with 6174 compute nodes, each with two 16-core Intel Haswell CPUs (Intel\textregistered Xeon\textregistered E5-2698 v3). The nodes are connected by a dragonfly network using the Aries interconnect.

% Conclusions
\section{Conclusions}\label{sec:conc}
The Fast Multipole Method, originally developed as a free-standing solver, can be effectively combined with Krylov iteration as a scalable and highly performant preconditioner for traditional low-order finite discretizations of elliptic boundary value problems.  In model problems it performs similarly to algebraic multigrid in convergence rate, while excelling in scalings where AMG becomes memory bandwidth-bound locally and/or syn-chronization-bound globally.  Additional algorithmic development and additional testing of implementations on emerging architectures are necessary to more fully define the niche in which FMM is the preconditioner of choice.  No preconditioner considered in isolation can address the fundamental architectural challenges of Kr-ylov methods for sparse linear systems, which are being simultaneously adapted to less synchronization tolerant computational environments through pipelining, but it is important to address the bottlenecks of preconditioning this most popular class of solvers by making a wide variety of tunable preconditioners available and better integrating them into the overall solver.  Fast multipole-based preconditioners are demonstrably ready to play an important role in the migration of sparse iterative solvers to the exascale.

%\newpage
\begin{acknowledgements}
The authors would like to acknowledge the open source software packages that made this work possible: PETSc~\cite{petsc-user-ref,petsc-web-page}, PetIGA~\cite{PetIGA} and IFISS~\cite{elman2007,ifiss}. We thank Dave Hewett, David May, Andy Wathen, and Ulrike Yang for helpful discussions  and comments and are indebted to Nathan Collier for his help with the PetIGA framework. We thank Lorena Barba for the support in developing ExaFMM. This publication was based on work supported in part by Award No KUK-C1-013-04, made by King Abdullah University of Science and Technology (KAUST). This work used the Extreme Science and Engineering Discovery Environment (XSEDE), which is supported by National Science Foundation grant number OCI-1053575.
\end{acknowledgements}

%%\bibliographystyle{spbasic}      % basic style, author-year citations
%\bibliographystyle{abbrv}
%%%\bibliography{./refs}

\bibliographystyle{abbrv}
\bibliography{./refs}

\begin{thebibliography}{10}

\bibitem{Aarseth1963}
S.~J. Aarseth.
\newblock Dynamical evolution of clusters of galaxies, {I}.
\newblock {\em Monthly Notices of the Royal Astronomical Society}, 126:233,
  1963.

\bibitem{adams2012}
M.~F. Adams.
\newblock A low memory, highly concurrent multigrid algorithm.
\newblock {\em arXiv:1207.6720v3}, 2012.

\bibitem{baker2012}
A.~H. Baker, R.~D. Falgout, T.~Gamblin, T.~V. Kolev, M.~Schulz, and U.~M. Yang.
\newblock Scaling algebraic multigrid solvers: On the road to exascale.
\newblock In C.~Bischof, H.-G. Hegering, W.~E. Nagel, and G.~Wittum, editors,
  {\em Competence in High Performance Computing}, pages 215--226. Springer,
  2012.

\bibitem{petsc-user-ref}
S.~Balay, J.~Brown, , K.~Buschelman, V.~Eijkhout, W.~D. Gropp, D.~Kaushik,
  M.~G. Knepley, L.~C. McInnes, B.~F. Smith, and H.~Zhang.
\newblock {PETS}c users manual.
\newblock Technical Report ANL-95/11 - Revision 3.4, Argonne National
  Laboratory, 2013.

\bibitem{petsc-web-page}
S.~Balay, J.~Brown, K.~Buschelman, W.~D. Gropp, D.~Kaushik, M.~G. Knepley,
  L.~C. McInnes, B.~F. Smith, and H.~Zhang.
\newblock {PETSc} {W}eb page, 2013.
\newblock http://www.mcs.anl.gov/petsc.

\bibitem{banerjee1979}
P.~K. Banerjee.
\newblock Non-linear problems of potential flow.
\newblock In {\em Developments in Boundary Element Methods---I}. Applied
  Science, 1979.

\bibitem{benzi2005}
M.~Benzi, G.~H. Golub, and J.~Liesen.
\newblock Numerical solution of saddle point problems.
\newblock {\em Acta Numerica}, 14:1--137, 2005.

\bibitem{borm2005}
S.~B{\"o}rm.
\newblock Hybrid cross approximation of integral operators.
\newblock {\em Numerische Mathematik}, 101:221--249, 2005.

\bibitem{brandt1977}
A.~Brandt.
\newblock Multi-level adaptive techniques ({MLAT}) for partial differential
  equations: Ideas and software.
\newblock In J.~R. Rice, editor, {\em Mathematical Software}, volume III, pages
  277--318. Academic Press, 1977.

\bibitem{brunton1996}
I.~Brunton.
\newblock {\em Solving Variable Coefficient Partial Differential Equations
  Using the Boundary Element Method}.
\newblock PhD thesis, University of Auckland, 1996.

\bibitem{cahouet1988}
J.~Cahouet and J.-P. Chabard.
\newblock Some fast 3{D} finite element solvers for the generalized {S}tokes
  problem.
\newblock {\em International Journal for Numerical Methods in Fluids},
  8:869--895, 1988.

\bibitem{carrier1988}
J.~Carrier, L.~Greengard, and V.~Rokhlin.
\newblock A fast adaptive multipole algorithm for particle simulations.
\newblock {\em SIAM Journal on Scientific and Statistical Computing},
  9(4):669--686, 1988.

\bibitem{chandrasekaran2006}
S.~Chandrasekaran, M.~Gu, and T.~Pals.
\newblock A fast {ULV} decomposition solver for hierarchically semiseparable
  representations.
\newblock {\em SIAM Journal on Matrix Analysis and Applications},
  28(3):603--622, 2006.

\bibitem{cheng1984}
A.~H.-D. Cheng.
\newblock Darcy's flow with variable permeability: A boundary integral
  solution.
\newblock {\em Water Resources Research}, 20:980--984, 1984.

\bibitem{clements1980}
D.~L. Clements.
\newblock A boundary integral equation method for the numerical solution of a
  second order elliptic equation with variable coefficients.
\newblock {\em Journal of the Australian Mathematical Society}, 22 (Series
  B):218--228, 1980.

\bibitem{concus}
P.~Concus and G.~H. Golub.
\newblock Use of fast direct methods for the efficient numerical solution of
  nonseparable elliptic equations.
\newblock {\em SIAM Journal on Numerical Analysis}, 10:1103--1120, 1973.

\bibitem{czechowski2012}
K.~Czechowski, C.~McClanahan, C.~Battaglino, K.~Iyer, P.-K. Yeung, and
  R.~Vuduc.
\newblock On the comunication complexity of 3d {FFT}s and its implications for
  exascale.
\newblock In {\em Proceedings of the 26th ACM International Conference on
  Supercomputing}, pages 205--214, 2012.

\bibitem{dehnen2002}
W.~Dehnen.
\newblock A hierarchical {O(N)} force calculation algorithm.
\newblock {\em Journal of Computational Physics}, 179(1):27--42, 2002.

\bibitem{demmel2008}
J.~Demmel, L.~Grigori, M.~F. Hoemmen, and J.~Langou.
\newblock Communication-avoiding parallel and sequential {QR} factorizations.
\newblock Technical Report EECS-2008-74, UC Berkeley, 2008.

\bibitem{dutt1996}
A.~Dutt, M.~Gu, and V.~Rokhlin.
\newblock Fast algorithms for polynomial interpolation integration and
  differentiation.
\newblock {\em SIAM Journal on Numerical Analysis}, 33(5):1689--1711, 1996.

\bibitem{elman2007}
H.~Elman, A.~Ramage, and D.~Silvester.
\newblock Algorithm {866}: {IFISS}, a {M}atlab toolbox for modelling
  incompressible flow.
\newblock {\em ACM Transactions on Mathematical Software}, 33:2--14, 2007.

\bibitem{elman2005}
H.~C. Elman, D.~J. Silvester, and A.~J. Wathen.
\newblock {\em Finite Elements and Fast Iterative Solvers: {W}ith applications
  in incompressible fluid dynamics}.
\newblock Oxford University Press, Oxford, 2005.

\bibitem{gahvari2011}
H.~Gahvari, A.~H. Baker, M.~Schulz, U.~M. Yang, K.~E. Jordan, and W.~Gropp.
\newblock Modeling the performance of an algebraic multigrid cycle on {HPC}
  platforms.
\newblock In {\em Proceedings of the international conference on
  Supercomputing}, pages 172--181, 2011.

\bibitem{Gholami2014}
A.~Gholami, D.~Malhotra, H.~Sundar, and G.~Biros.
\newblock {FFT}, {FMM}, or {MULTIGRID}? {A} comparative study of
  state-of-the-art poisson solvers.
\newblock {\em arXiv:1408.6497}, 2014.

\bibitem{golub1996}
G.~H. Golub and C.~F. van Loan.
\newblock {\em Matrix Computations}.
\newblock The John Hopkins University Press, Baltimore, 1996.

\bibitem{golub1961}
G.~H. Golub and R.~S. Varga.
\newblock Chebyshev semi iterative methods, successive overrelaxation iterative
  methods and second order {R}ichardson iterative methods. {P}art {I}.
\newblock {\em Numerische Mathematik}, 3:147--156, 1961.

\bibitem{graham2007}
I.~G. Graham, P.~O. Lechner, and R.~Scheichl.
\newblock Domain decomposition for multiscale {PDE}s.
\newblock {\em Numerische Mathematik}, 106:589--626, 2007.

\bibitem{grasedyck2003}
L.~Grasedyck and W.~Hackbusch.
\newblock Construction and arithmetics of {H}-matrices.
\newblock {\em Computing}, 70:295--334, 2003.

\bibitem{grasedyck2008}
L.~Grasedyck, W.~Hackbusch, and R.~Kriemann.
\newblock Performance of {H-LU} preconditioning for sparse matrices.
\newblock {\em Computational Methods in Applied Mathematics}, 8(4):336--349,
  2008.

\bibitem{greenbaum1997}
A.~Greenbaum.
\newblock {\em Iterative Methods for Solving Linear Systems}.
\newblock SIAM, USA, 1997.

\bibitem{greengard2009}
L.~Greengard, D.~Gueyffier, P.~G. Martinsson, and V.~Rokhlin.
\newblock Fast direct solvers for integral equations in complex three
  dimensional domains.
\newblock {\em Acta Numerica}, 18:243--275, 2009.

\bibitem{greengard1987}
L.~Greengard and V.~Rokhlin.
\newblock A fast algorithm for particle simulations.
\newblock {\em Journal of Computational Physics}, 73(2):325--348, 1987.

\bibitem{gu1996}
M.~Gu and S.~C. Eisenstat.
\newblock Efficient algorithms for computing a strong rank-revealing qr
  factorization.
\newblock {\em SIAM Journal on Scientific Computing}, 17(4):848--869, 1996.

\bibitem{gumerov2008}
N.~A. Gumerov and R.~Duraiswami.
\newblock Fast multipole methods on graphics processors.
\newblock {\em Journal of Computational Physics}, 227:8290--8313, 2008.

\bibitem{hackbusch1999}
W.~Hackbusch.
\newblock A sparse matrix arithmetic based on {H}-matrices, {P}art {I}:
  Introduction to {H}-matrices.
\newblock {\em Computing}, 62:89--108, 1999.

\bibitem{hackbusch2000onh2}
W.~Hackbusch, B.~Khoromskij, and S.~Sauter.
\newblock On {H}2-{M}atrices.
\newblock In H.-J. Bungartz, R.~Hoppe, and C.~Zenger, editors, {\em Lectures on
  Applied Mathematics}, pages 9--29. Springer Berlin Heidelberg, 2000.

\bibitem{henson2002}
V.~E. Henson and U.~M. Yang.
\newblock Boomer{AMG}: A parallel algebraic multigrid solver and
  preconditioner.
\newblock {\em Applied Numerical Mathematics}, 41:155--177, 2002.

\bibitem{hestenes1952}
M.~R. Hestenes and E.~Stiefel.
\newblock Methods of conjugate gradients for solving linear systems.
\newblock {\em Journal of Research of the National Bureau of Standards},
  49:409--436, 1952.

\bibitem{ho2012}
K.~L. Ho and L.~Greengard.
\newblock A fast direct solver for structured linear systems by recursive
  skeletonization.
\newblock {\em SIAM Journal on Scientific Computing}, 34(5):A2507--A2532, 2012.

\bibitem{hundsdorfer2003}
W.~Hundsdorfer and J.~G. Verwer.
\newblock {\em Numerical Solution of Time-Dependent
  Advection-Diffusion-Reaction Equations}.
\newblock Springer-Verlag, Germany, 2003.

\bibitem{kahan1965}
W.~Kahan.
\newblock Pracniques: Further remarks on reducing truncation errors.
\newblock {\em Communications of the ACM}, 8(1):40, 1965.

\bibitem{kong2011}
W.~Y. Kong, J.~Bremer, and V.~Rokhlin.
\newblock An adaptive fast direct solver for boundary integral equations in two
  dimensions.
\newblock {\em Applied and Computational Harmonic Analysis}, 31:346--369, 2011.

\bibitem{langer2007}
U.~Langer, G.~Of, O.~Steinbach, and W.~Zulehner.
\newblock Inexact fast multipole boundary element tearing and interconnecting
  methods.
\newblock In {\em Domain Decomposition Methods in Science and Engineering XVI}.
  Springer Berlin Heidelberg, 2007.

\bibitem{lashuk2012}
I.~Lashuk, A.~Chandramowlishwaran, H.~Langston, T.-A. Nguyen, R.~Sampath,
  A.~Shringarpure, R.~Vuduc, L.~Ying, D.~Zorin, and G.~Biros.
\newblock A massively parallel adaptive fast multipole method on heterogeneous
  architectures.
\newblock {\em Communications of the ACM}, 55(5):101--109, 2012.

\bibitem{liberty2007}
E.~Liberty, F.~Woolfe, P.~G. Martinsson, V.~Rokhlin, and M.~Tygert.
\newblock Randomized algorithms for the low-rank approximation of matrices.
\newblock {\em PNAS}, 104(51):20167--20172, 2007.

\bibitem{dhairya2015}
D.~Malhotra and G.~Biros.
\newblock {PVFMM}: A parallel kernel independent {FMM} for particle and volume
  potentials.
\newblock {\em Communications in Computational Physics}, 18:808--830, 9 2015.

\bibitem{martin2004}
R.~M. Martin.
\newblock {\em Electronic Structure: Basic Theory and Practical Methods}.
\newblock Cambridge University Press, Cambridge, UK, 2004.

\bibitem{meijerink1977}
J.~A. Meijerink and H.~A. van~der Vorst.
\newblock An iterative solution method for linear systems of which the
  coefficient matrix is a symmetric {$M$}-matrix.
\newblock {\em Mathematics of Computation}, 31:148--162, 1977.

\bibitem{PetIGA}
V.~C. N.~Collier, L.~Dalcin.
\newblock {PetIGA}: High-performance isogeometric analysis.
\newblock {\em arxiv}, (1305.4452), 2013.
\newblock http://arxiv.org/abs/1305.4452.

\bibitem{paige1975}
C.~C. Paige and M.~A. Saunders.
\newblock Solution of sparse indefinite systems of linear equations.
\newblock {\em SIAM Journal on Numerical Analysis}, 12:617--629, 1975.

\bibitem{pan2000}
C.~T. Pan.
\newblock On the existence and computation of rank-revealing {LU}
  factorizations.
\newblock {\em Linear Algebra and its Applications}, 316:199--222, 2000.

\bibitem{rjasanow2002}
S.~Rjasanow.
\newblock Adaptive cross approximation of dense matrices.
\newblock In {\em International Association for Boundary Element Methods}, UT
  Austin, TX, USA, May 28-30 2002.

\bibitem{saad2003}
Y.~Saad.
\newblock {\em Iterative Methods for Sparse Linear Systems}.
\newblock SIAM, Philadelphia, 2003.

\bibitem{saad1986}
Y.~Saad and M.~H. Schultz.
\newblock {GMRES}: {A} generalized minimal residual algorithm for solving
  nonsymmetric linear systems.
\newblock {\em SIAM Journal on Scientific and Statistical Computing},
  7:856--869, 1986.

\bibitem{sambavaram2003}
S.~R. Sambavaram, V.~Sarin, A.~Sameh, and A.~Grama.
\newblock Multipole-based preconditioners for large sparse linear systems.
\newblock {\em Parallel Computing}, 29:1261--1273, 2003.

\bibitem{sauter2011}
S.~A. Sauter and {Ch.~Schwab}.
\newblock {\em Boundary Element Methods}.
\newblock Springer-Verlag, Heidelberg, 2011.

\bibitem{ifiss}
D.~Silvester, H.~Elman, and A.~Ramage.
\newblock {I}ncompressible {F}low and {I}terative {S}olver {S}oftware ({IFISS})
  version 3.2, May 2012.
\newblock \\ {\tt http://www.manchester.ac.uk/ifiss}.

\bibitem{Teng1998}
S.-H. Teng.
\newblock Provably good partitioning and load balancing algorithms for parallel
  adaptive {N}-body simulation.
\newblock {\em SIAM Journal on Scientific Computing}, 19(2):635--656, 1998.

\bibitem{trottenberg2001}
U.~Trottenberg, C.~Oosterlee, and A.~Sch\"{u}ller.
\newblock {\em Multigrid}.
\newblock Academic Press, London, 2001.

\bibitem{wardle1978}
L.~J. Wardle and J.~M. Crotty.
\newblock Two dimensional boundary integral equation analysis for
  non-homogeneous mining applications.
\newblock In {\em Recent Advances in Boundary Element Methods}. Pentch Press,
  1978.

\bibitem{warren1992}
M.~S. Warren and J.~K. Salmon.
\newblock Astrophysical {N}-body simulation using hierarchical tree data
  structures.
\newblock In {\em Proceedings of the 1992 ACM/IEEE Conference on
  Supercomputing}, pages 570--576, 1992.

\bibitem{warren1995}
M.~S. Warren and J.~K. Salmon.
\newblock A portable parallel particle program.
\newblock {\em Computer Physics Communications}, 87:266--290, 1995.

\bibitem{wathen2009}
A.~Wathen and T.~Rees.
\newblock Chebyshev semi-iteration in preconditioning for problems including
  the mass matrix.
\newblock {\em Electronic Transactions on Numerical Analysis}, 34:125--135,
  2009.

\bibitem{wathen1987}
A.~J. Wathen.
\newblock Realistic eigenvalue bounds for the {G}alerkin mass matrix.
\newblock {\em IMA Journal of Numerical Analysis}, 7:449--457, 1987.

\bibitem{yokota2013a}
R.~Yokota.
\newblock An {FMM} based on dual tree traversal for many-core architectures.
\newblock {\em Journal of Algorithms and Computational Technology},
  7(3):301--324, 2013.

\end{thebibliography}

%\begin{thebibliography}{}
%\end{thebibliography}

\end{document}